\documentclass[a4paper]{cedram-smai-jcm}

\usepackage{inputenc}
\usepackage{graphicx}
\usepackage{subfigure}
\usepackage{caption}
\usepackage{pgf,tikz}
\usepackage{mathrsfs}
\usepackage{color}
\usetikzlibrary{arrows,decorations.markings,shapes.misc,decorations.pathreplacing,shapes.geometric,shapes.arrows}
\usepackage[english]{babel}
\usepackage{amsmath,amssymb,bbm,xcolor,enumerate}
\usepackage{array}
\usepackage[lined,boxed,commentsnumbered]{algorithm2e}

\tikzset{
  on each segment/.style={
    decorate,
    decoration={
      show path construction,
      moveto code={},
      lineto code={
        \path [#1]
        (\tikzinputsegmentfirst) -- (\tikzinputsegmentlast);
      },
      curveto code={
        \path [#1] (\tikzinputsegmentfirst)
        .. controls
        (\tikzinputsegmentsupporta) and (\tikzinputsegmentsupportb)
        ..
        (\tikzinputsegmentlast);
      },
      closepath code={
        \path [#1]
        (\tikzinputsegmentfirst) -- (\tikzinputsegmentlast);
      },
    },
  },
  mid arrow/.style={postaction={decorate,decoration={
        markings,
        mark=at position .5 with {\arrow[#1,thick]{stealth}}
        }}},
}

\catcode`\|=\active \def|{
\fontencoding{T1}\selectfont\symbol{124}\fontencoding{\encodingdefault}}

\newcommand{\tmdate}[1]{\today}

\newcommand{\assign}{:=}

\reversemarginpar

\title[Flip procedure in approximation of multiple-component shapes]{Flip procedure in geometric approximation of multiple-component shapes -- Application to multiple-inclusion detection}

\author[P. Bonnelie]{\firstname{Pierre} \lastname{Bonnelie}}
\address{Institut de recherche XLIM. P\^ole Math\'ematiques-Informatique-Image. UMR CNRS 7252. Universit\'e de Limoges, France.}
\email{pierre.bonnelie@etu.unilim.fr}
\thanks{The first author is granted by Labex SIGMA-LIM}

\author[L. Bourdin]{\firstname{Lo\"ic} \lastname{Bourdin}}
\address{Institut de recherche XLIM. P\^ole Math\'ematiques-Informatique-Image. UMR CNRS 7252. Universit\'e de Limoges, France.}
\email{loic.bourdin@unilim.fr}

\author[F. Caubet]{\firstname{Fabien} \lastname{Caubet}}
\address{Institut de Math\'ematiques de Toulouse. UMR CNRS 5219. Universit\'e de Toulouse, France.}
\email{fabien.caubet@math.univ-toulouse.fr}

\author[O. Ruatta]{\firstname{Olivier} \lastname{Ruatta}}
\address{Institut de recherche XLIM. P\^ole Math\'ematiques-Informatique-Image. UMR CNRS 7252. Universit\'e de Limoges, France.}
\email{olivier.ruatta@unilim.fr}

  

\def\R{\mathbb{R}}
\def\N{\mathbb{N}}
\def\N{\mathbb{N}}

\newcommand{\argmin}{{\rm argmin}\,}

\newcommand{\nn}{\boldsymbol{\mathrm{n}}}
\newcommand{\DD}{{\rm D}}
\newcommand{\HH}{{\rm H}}

\newcommand{\LL}{{\rm L}}

\newcommand{\GII}{\boldsymbol{\mathrm{I}}}

\newcommand{\GWW}{\boldsymbol{\mathrm{W}}}
\newcommand{\Gn}{\boldsymbol{\mathrm{n}}}

\newcommand{\GU}{\boldsymbol{U}}
\newcommand{\GV}{\boldsymbol{V}}

\newcommand{\0}{\boldsymbol{0}}
\newcommand{\Gtheta}{\boldsymbol{\mathrm{\theta}}}

\newcommand{\abs}[1]{\left\vert #1 \right\vert}
\newcommand{\priv}[2]{#1  \backslash  \overline{#2}}

\newcommand{\BE}{\begin{equation}}
\newcommand{\EE}{\end{equation}}
\newcommand{\BA}{\begin{array}}
\newcommand{\EA}{\end{array}}
\newcommand{\BEAN}{\begin{eqnarray*}}
\newcommand{\EEAN}{\end{eqnarray*}}
\newcommand{\BM}{\begin{multline}}
\newcommand{\EM}{\end{multline}}

\newcommand{\Oomega}{\priv{\Omega}{\omega}}
\newcommand{\Oomegat}{\priv{\Omega}{\omega_{t}}}

\renewcommand{\(}{\left(}
\renewcommand{\)}{\right)}
\renewcommand{\[}{\left[}

\def \restriction#1#2{\mathchoice
              {\setbox1\hbox{${\displaystyle #1}_{\scriptstyle #2}$}
              \restrictionaux{#1}{#2}}
              {\setbox1\hbox{${\textstyle #1}_{\scriptstyle #2}$}
              \restrictionaux{#1}{#2}}
              {\setbox1\hbox{${\scriptstyle #1}_{\scriptscriptstyle #2}$}
              \restrictionaux{#1}{#2}}
              {\setbox1\hbox{${\scriptscriptstyle #1}_{\scriptscriptstyle #2}$}
              \restrictionaux{#1}{#2}}}
\def\restrictionaux#1#2{{#1\,\smash{\vrule height .8\ht1 depth .85\dp1}}_{\,#2}}

\begin{document}

\begin{abstract}
We are interested in geometric approximation by parameterization of two-dimensional multiple-component shapes, in particular when the number of components is \textit{a priori} unknown. 
Starting a standard method based on successive shape deformations with a one-component initial shape in order to approximate a multiple-component target shape usually leads the deformation flow to make the boundary evolve until it surrounds all the components of the target shape. This classical phenomenon tends to create double points on the boundary of the approximated shape. 

In order to improve the approximation of multiple-component shapes (without any knowledge on the number of components in advance), we use in this paper a piecewise B\'ezier parameterization and we consider two procedures called \textit{intersecting control polygons detection} and \textit{flip procedure}. The first one allows to prevent potential collisions between two parts of the boundary of the approximated shape, and the second one permits to change its topology by dividing a one-component shape into a two-component shape. 


For an experimental purpose, we include these two processes in a basic geometrical shape optimization algorithm and test it on the classical inverse obstacle problem. This new approach allows to obtain a numerical approximation of the unknown inclusion, detecting both the topology (\textit{i.e.} the number of connected components) and the shape of the obstacle. Several numerical simulations are performed.
\end{abstract}

\maketitle

\section{Introduction}

Geometric shape approximation methods are commonly based on successive shape deformations, where the boundary of the approximated shape is parameterized and evolves at each step in a direction given by the deformation flow. This technique is widely used for example in shape optimization problems where the flow is given by the so-called \textit{shape gradient} (see, e.g., Chapter 5 of the book \cite{HP} of Henrot \textit{et al.}), or in image segmentation (see, e.g., \cite{Kichenassamy94gradientflows}). Numerous parameterizations of the boundary have been considered in the literature, such as polygons, Fourier series, etc. Each of these parameterizations has its own advantages and drawbacks, that depend on the nature of the problem studied.

\begin{figure}[h]
\subfigure[Two-dimensional case]{
\begin{tikzpicture}[scale=0.4, every node/.style={scale=0.4}]
\draw[red] (2,0) .. controls (1.5,-0.3) and (-0.7,0.2) .. (-1,1);
\draw[red] (-1,1) .. controls (-1.7,1.9) and (-2.5,2.3) .. (-3,2);
\draw[red] (-3,2) .. controls (-5.5,1) and (-4.5,-4) .. (-1.5,-2.2);
\draw[red] (-1.5,-2.2) .. controls (1,-0.4) and (1,-0.2) .. (2,-1.2);
\draw[red] (2,-1.2) .. controls (6,-4) and (7,3) .. (2,0);
\draw[blue,thick] (4,-0.5) circle (1.2cm);
\draw[blue,thick] (-2.5,0) circle (1.6cm);
\draw[blue,thick] (0,2.7) -- (2,2.7);
\node (target) at (3.8,2.7) {\LARGE Target shape};
\draw[red] (0,1.7) -- (2,1.7);
\node (final) at (4.8,1.7) {\LARGE Approximated shape};
\end{tikzpicture}
}
\subfigure[Three-dimensional case]{
\includegraphics[scale=0.1]{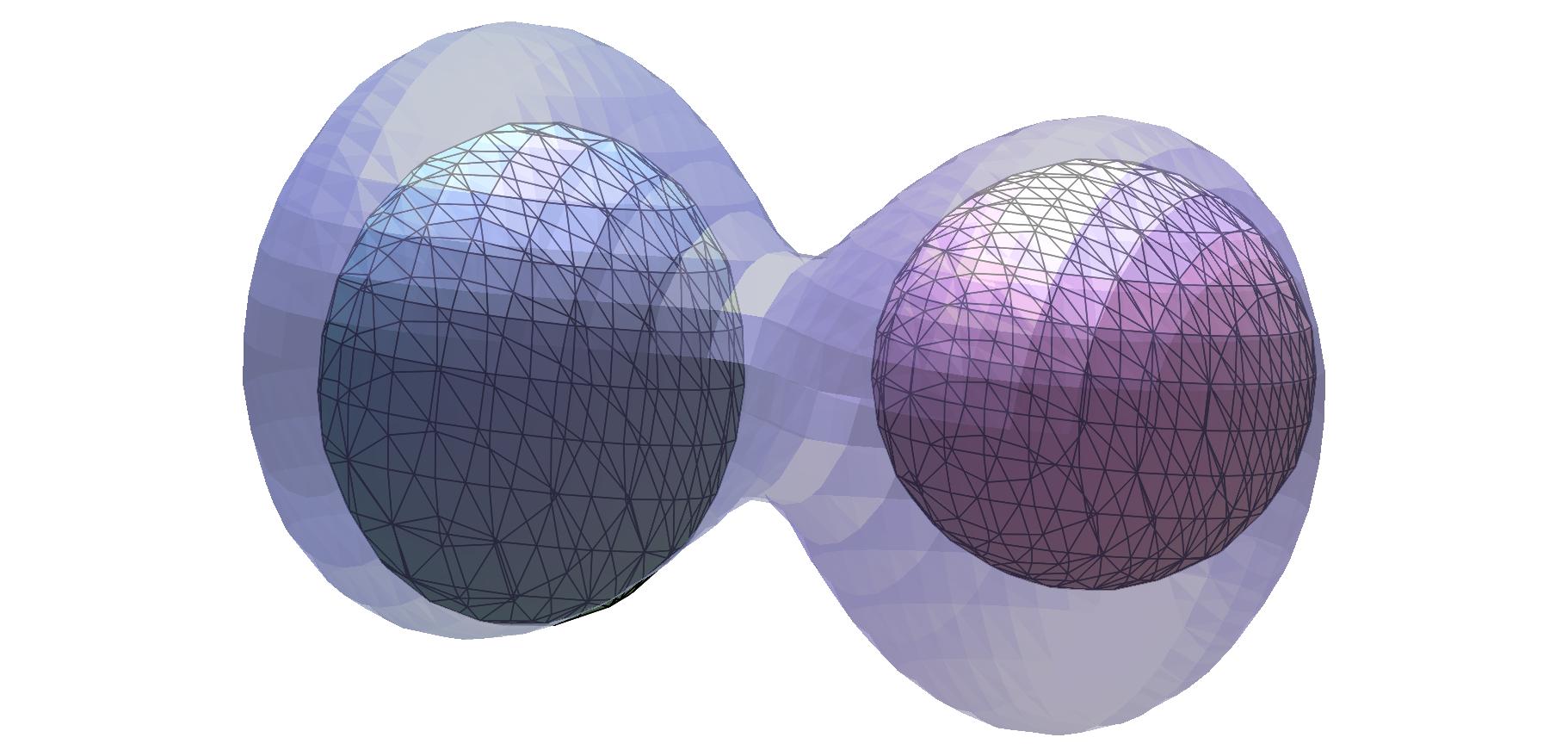}
}
\caption{Geometric shape approximation of a two-component target shape starting from a one-component initial shape.}
\label{figIntro}
\end{figure}
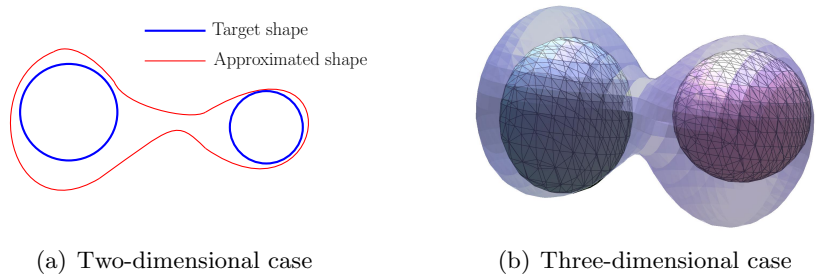

In this paper we are especially interested in the geometric approximation of multiple-component shapes, in particular when the number of components is \textit{a priori} unknown. Starting a parameterization method with a one-component initial shape in order to approximate a multiple-component target shape usually leads the deformation flow to make the boundary evolve until it surrounds all the components of the target shape (see Figure~\ref{figIntro} for illustrations). This classical phenomenon tends to create double points on the boundary of the approximated shape. 

In order to improve the approximation of multiple-component shapes, our idea is to look for an appropriate parameterization that allows to achieve two numerical tasks. Firstly the parameterization has to be well-suited in order to prevent the potential formation of double points, \textit{i.e.} to locate the parts of the boundary that are close to each other. Secondly it has to be adapted in order to easily change the topology of the approximated shape, precisely in order to divide a one-component shape into a two-component shape. Moreover, for practical uses, we look for a complete method that is easily implementable with a relatively low numerical cost.

We present in this paper a method based on a B\'ezier parameterization. The main idea is that this polynomial parameterization can be approximated by its control polygon. In particular one can easily prevent the potential formation of double points by looking for intersecting control polygons. We refer to Section~\ref{intersectdetection} for details on the so-called \textit{intersecting control polygons detection}. Once this first step is achieved, one can easily reorganize the control points of the B\'ezier parameterization in order to modify the topology of the shape, precisely in order to divide one component into two. We refer to Section~\ref{flipalgo} for details on the so-called \textit{flip procedure}.\footnote{Actually a similar procedure can also be considered in order to merge two components into one (see Appendix~\ref{2for1} for some details).} In this work we detail the above method in the two-dimensional case, using piecewise B\'ezier curves.\footnote{Let us give a brief discussion on the three-dimensional case. We refer for instance to~\cite{MsPhan14} where deformation of piecewize B\'ezier surfaces is presented with an implementation. Note that the adaptation of the complete algorithmic setting of the flip procedure to the three-dimensional case would be nontrivial since it would increase the algorithmic and combinatoric complexities. Numerous considerations about this generalization could be addressed, however we postpone this interesting issue to a future work.}



\begin{figure}[h]
\centering
\begin{tikzpicture}[scale=0.4, every node/.style={scale=0.4}]
\draw (1.,4) circle (7.cm);
\draw[green] (1.5,6) circle (2.7cm);
\draw[blue] (2,6) .. controls (0,5.9) and (-0.5,5.) .. (-1,5);
\draw[blue] (-1,5) .. controls (-1.7,4.9) and (-2.5,5.3) .. (-3,5);
\draw[blue] (-3,5) .. controls (-5.5,4) and (-4.5,-1) .. (-1.5,0.2);
\draw[blue] (-1.5,0.2) .. controls (1,1.7) and (1.5,1.) .. (2,1);
\draw[blue] (2,1) .. controls (5,1) and (7,6) .. (2,6);
\draw[red] (2.7,6) .. controls (0,5.9) and (-0.,4.3) .. (-1.2,4.8) ;
\draw[red] (-1.2,4.8) .. controls (-1.7,4.9) and (-2.5,5.3) .. (-3.5,5);
\draw[red] (-3.5,5) .. controls (-5.5,4) and (-4.5,-1) .. (-1.3,1.2);
\draw[red] (-1.3,1.2) .. controls (1,2.7) and (1.5,0.2) .. (2.3,0.4);
\draw[red] (2.3,0.4) .. controls (5.9,1) and (7,6) .. (2.7,6);
\draw (9,4.8) -- (11,4.8);
\node (target) at (13.95,4.8) {\LARGE Exterior boundary $\partial \Omega$};
\draw[blue] (9,4) -- (11,4);
\node (target) at (15.4,4) {\LARGE Target shape (or exact shape) $\omega_{\rm ex}$};
\draw[green] (9,3.2) -- (11,3.2);
\node (init) at (14.6,3.2) {\LARGE Initial approximated shape};
\draw[red] (9,2.4) -- (11,2.4);
\node (final) at (14.5,2.4) {\LARGE Final approximated shape};
\end{tikzpicture}
\caption{Illustration of reconstruction for the inverse obstacle problem.}
\label{FigReconstructionIntro}
\end{figure}
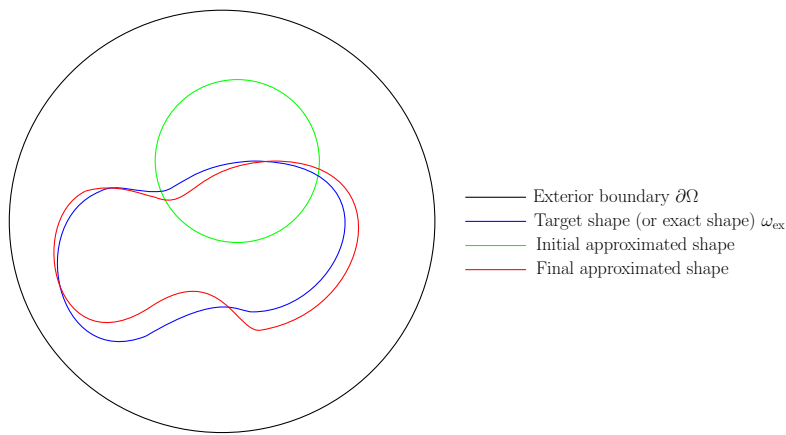

In order to test the two procedures introduced in this paper, we perform numerical simulations on the classical inverse obstacle problem. Precisely we consider the inverse problem of detecting some unknown inclusions~$\omega_{\rm ex}$ in a larger bounded domain $\Omega$ from boundary measurements made on $\partial\Omega$. The aim is to reconstruct numerically an approximation of the target shape~$\omega_{\rm ex}$ using shape optimization tools (see Figure~\ref{FigReconstructionIntro} for an illustration). In this paper we will study this inverse problem by minimizing a shape least-square functional.

We briefly recall now the major shape optimization techniques used in order to study this problem in the literature. 
Two main categories are {\it topological} and {\it geometric} shape optimization methods. The \textit{topological} gradient approach was introduced by Schumacher in~\cite{Sch95} and Sokolowski~\textit{et al.} in~\cite{SokZoc99}. 
This method is based on asymptotic expansions and consequently is essentially adapted for relatively small inclusions. Moreover, even if the topological optimization is useful in order to find the number of inclusions, it may be not well-suited in order to find a satisfactory approximation of the shape of the inclusions (see, e.g.,~\cite{Caubet12} and references therein). { We refer to \cite{AmmKan04} and references therein for a comprehensive mathematical treatment with theoretical and numerical results about reconstruction of small inclusions from boundary measurements}. We also refer to the work~\cite{GarGui01} for generalities on topological asymptotic expansion in the elastic context. 
In the \textit{geometric} shape optimization category, two main techniques are addressed in the literature. They are both based on the computation of a shape gradient used as a flow making the shape evolve. These two methods use different representations of the shape and different techniques to deform it. 
The first approach is the so-called level set approach (see, e.g., the survey~\cite{BurOsh05} of Burger~\textit{et al.} and references therein or~\cite{Set99}). It is originally based on an implicit representation of the approximated shape on a fixed mesh and, in the case of inverse problems, some regularization methods are usually needed (as curve shortening in, e.g.,~\cite[Section~8]{San96}). In order to detect several inclusions, this method does not need any \textit{a priori} knowledge on the number of inclusions.
The second approach is based on { boundary variations {\it via} mesh variations} and, in the case of inverse problems, on an explicit representation of the approximated shape. This method is used e.g. in the work~\cite{Dambrine} of Afraites~\textit{et al.} (where a regularization by parameterization is used). Note that the standard algorithm based on shape derivatives moving the mesh does not provide the opportunity to change the topology of the shape and consequently the number of inclusions has to be known in advance.
Recent works propose to mix several of the above different approaches. For instance we refer to the works of Allaire~\textit{et al.} in~\cite{AllDap14,AllGou05} and Burger~\textit{et al.} in~\cite{BurHac04} (see also the thesis~\cite[Section~5]{Persson:2005:MGI:1087686}) that combine the classical geometric shape optimization through the level set method and the topological gradient, to the work of Pantz~\textit{et al.} in~\cite{PanTra07} which develops an algorithm using boundary variations, topological derivatives and homogenization methods and to the works of Caubet~\textit{et al.} in~\cite{caubet:hal-01191099} and Christiansen~\textit{et al.} in~\cite{Christiansen2014} which couple topological and boundary variations approaches. 

The method presented in this paper is based only on mesh variation techniques. The parameterization by piecewise B\'ezier curves and the flip procedure permit to dynamically change the topology of the shape in order to find the number of inclusions, and the shape derivatives approach allows to approximate the shape of the inclusions with an explicit representation. This new method seems to be well-suited in order to study the above inverse obstacle problem, in particular in the case where the number of inclusions is \textit{a priori} unknown, and can be seen as an alternative to the above mixed methods combining the level set approach and the topological gradient.


\paragraph{Organization of the paper.}
The paper is organized as follows. 
Section~\ref{SectionBezier} recalls some basics and notations about piecewise B\'ezier curves. Section~\ref{SectionFlip} is concerned with the two main features of this paper, that is, the intersecting control polygons detection and the flip procedure. 
Section~\ref{SectionObjectsDetection} is dedicated to several numerical simulations in the context of the inverse obstacle problem. 

\section{Notations and basics on piecewise B\'ezier curves} \label{SectionBezier}
\begin{figure}[h]
\definecolor{mycolor}{RGB}{240,240,240}
\centering
\begin{tikzpicture}[scale=0.6, every node/.style={scale=0.6}]
\draw[dashed,fill=mycolor] (0,-1) -- (1,2) -- (3,3) -- (6,2) -- (7,0) -- cycle;
\draw [red,  domain=0:1, samples=40]
 plot ({0*(1-\x)^4 + 1*4*\x*(1-\x)^3 + 3*6*\x^2*(1-\x)^2 + 6*4*\x^3*(1-\x) + 7*\x^4}, {(-1)*(1-\x)^4 + 2*4*\x*(1-\x)^3 + 3*6*\x^2*(1-\x)^2 + 2*4*\x^3*(1-\x) + 0*\x^4} );
\draw[dashed] (0,-1) -- (1,2) -- (3,3) -- (6,2) -- (7,0);
\node[circle,fill=black,draw=black,scale=0.3] (p0) at (0,-1) {};
\node at (-0.5,-1) {$P_{0}$};
\node[circle,fill=black,draw=black,scale=0.3] (p1) at (1,2) {};
\node at (0.5,2.2) {$P_{1}$};
\node[circle,fill=black,draw=black,scale=0.3] (p2) at (3,3) {};
\node at (3,3.4) {$P_{2}$};
\node[circle,fill=black,draw=black,scale=0.3] (p3) at (6,2) {};
\node at (6.5,2.2) {$P_{3}$};
\node[circle,fill=black,draw=black,scale=0.3] (p4) at (7,0) {};
\node at (7.5,0) {$P_{4}$};
\end{tikzpicture}
\caption{A non-closed B\'ezier curve of degree $4$ lying in the convex hull of its control polygon.}
\label{beziercurve}
\end{figure}
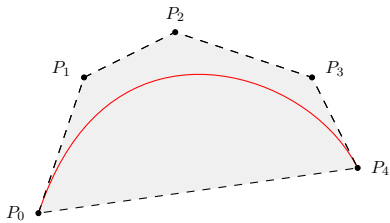

In this section we fix our notations and recall some basics about B\'ezier curves (see, e.g., \cite{Farin,SederbergCourseNotes} or \cite[from p.~409]{frey2008} for more details). Let $d \in \mathbb{N}^*$ and a set of $d+1$ points $P_0,\ldots,P_d$ of $\mathbb{R}^2$. The associated B\'ezier curve, denoted by $B([P_0,\ldots,P_d])$, is defined by
$$ \forall t \in [0,1], \quad B([P_0,\ldots,P_d],t) := \displaystyle \sum_{j=0}^{d} P_{j} b_{j,d}(t) , $$
where $b_{j,d}$ are the classical Bernstein polynomials given by
\begin{equation*} \label{Bernstein}
 b_{j,d}(t) := \binom{d}{j} t^j (1-t)^{d-j}.
\end{equation*}
The integer $d$ is the \textit{degree} of the curve and the points $P_0, \ldots, P_d$ are its \textit{control points} (or its \textit{control polygon}). Note that a B\'ezier curve does not go through its control points in general. However it starts at $P_0$ and finishes at $P_d$. If $P_0=P_d$, the B\'ezier curve is said to be \textit{closed}. Each point of a B\'ezier curve is a convex combination of its control points. As a consequence, a B\'ezier curve lies in the convex hull of its control polygon (see Figure~\ref{beziercurve}). 

\begin{remark}
As B\'ezier curves are widely used in Computer Aided Geometric Design (see \cite{Farin,SederbergCourseNotes}), they are commonly defined as parametric curves lying in the euclidean space $\mathbb{R}^2$ (or $\mathbb{R}^3$). However this definition can be extended to $\mathbb{R}^n$ for any $n \in \mathbb{N}^*$. In this paper, we are only interested in the two-dimensional case $n=2$.
\end{remark}

In this paper we focus on the geometric approximation of boundaries of two-dimensional bounded shapes with the help of B\'ezier curves. In the sequel no distinction will be done between a two-dimensional bounded shape and its boundary. 

Using a single closed B\'ezier curve in order to approximate a two-dimensional shape is not an efficient method for several reasons. Indeed, in order to approximate a shape with a lot of geometric features, one would need to increase the number of degrees of freedom, \textit{i.e.} the number of control points. However, as is very well-known, increasing the degree of an approximating polynomial curve leads to a classical oscillation phenomenon and, in the particular case of a B\'ezier polynomial curve, it leads to numerical instabilities (due to the ill-conditionness of the Bernstein-Vandermonde matrices, see, e.g.,~\cite{MARCO2007616}). Moreover, since each control point has a global influence on the curve, one could not handle local complexities of a shape with a single B\'ezier curve. The classical idea is then to divide the curve in several B\'ezier curves of small degrees. This leads us to recall the following definition of piecewise B\'ezier curves.


Let $N \in \mathbb{N}^*$, $d \in \mathbb{N}^*$ and a set of $N(d+1)$ control points $P_{1,0},\ldots, P_{1,d},\ldots, P_{N,d}$ of $\mathbb{R}^2$ satisfying the \textit{continuity relations} $P_{i,d} = P_{i+1,0}$ for every $i=1,\ldots,N-1$.\footnote{The continuity relations guarantee the well-definedness and the continuity of the piecewise B\'ezier curve.} The associated piecewise B\'ezier curve, denoted by $B([P_{1,0},\ldots,P_{N,d}])$\footnote{One would note here a conflict in notations of a B\'ezier curve and of a piecewise B\'ezier curve. In the sequel no confusion is possible since we will only consider piecewise B\'ezier curves.}, is defined by
\begin{multline*}
\label{defpbc}
\forall t \in [0,1], \quad B([P_{1,0},\ldots,P_{N,d}],t) := B([P_{i,0},\ldots,P_{i,d}],Nt-i+1), \\ \textrm{ if } t \in \left[ \frac{i-1}{N} , \frac{i}{N} \right],
\textrm{ }i \textrm{ ranges from } 1 \textrm{ to } N.
\end{multline*}
The global curve is then composed of $N$ B\'ezier curves called \textit{patches}. 
Note that a piecewise B\'ezier curve goes through $P_{i,0}$ and $P_{i,d}$ for all $i=1,\ldots,N$. If $P_{1,0}=P_{N,d}$, the piecewise B\'ezier curve is said to be \textit{closed}. 


\begin{remark}\label{remarkd3}
In practice we use cubic patches ($d=3$) because they are sufficient in order to recover many geometrical situations, such as inflexion points (see Figure~\ref{piecewiseBeziercurve}).
\begin{figure}[h]
\begin{tikzpicture}[scale=0.6, every node/.style={scale=0.6}]
\node[circle,fill=black,draw=black,scale=0.3] (p0) at (0,4) {};
\node[circle,fill=black,draw=black,scale=0.3] (p1) at (0,3) {};
\node[circle,fill=black,draw=black,scale=0.3] (p2) at (1,2) {};
\node[circle,fill=black,draw=black,scale=0.3] (p3) at (2,3) {};
\node[circle,fill=black,draw=black,scale=0.3] (p4) at (3,4) {};
\node[circle,fill=black,draw=black,scale=0.3] (p5) at (4,3) {};
\node[circle,fill=black,draw=black,scale=0.3] (p6) at (4,2) {};
\node[circle,fill=black,draw=black,scale=0.3] (p7) at (5,0) {};
\node[circle,fill=black,draw=black,scale=0.3] (p8) at (7,0) {};
\node[circle,fill=black,draw=black,scale=0.3] (p9) at (8,1) {};
\node[circle,fill=black,draw=black,scale=0.3] (p10) at (9,2) {};
\node[circle,fill=black,draw=black,scale=0.3] (p11) at (10,4) {};
\node[circle,fill=black,draw=black,scale=0.3] (p12) at (8,4) {};
\node[circle,fill=black,draw=black,scale=0.3] (p13) at (6,4) {};
\node[circle,fill=black,draw=black,scale=0.3] (p14) at (8,6) {};
\node[circle,fill=black,draw=black,scale=0.3] (p15) at (6,6) {};
\node[circle,fill=black,draw=black,scale=0.3] (p16) at (5,7) {};
\node[circle,fill=black,draw=black,scale=0.3] (p17) at (3,7) {};
\node[circle,fill=black,draw=black,scale=0.3] (p18) at (2,6) {};
\node[circle,fill=black,draw=black,scale=0.3] (p19) at (1,6) {};
\node[circle,fill=black,draw=black,scale=0.3] (p20) at (0,5) {};
\draw[thick] (p0) .. controls (p1) and (p2) .. (p3);
\draw[thick] (p3) .. controls (p4) and (p5) .. (p6);
\draw[thick] (p6) .. controls (p7) and (p8) .. (p9);
\draw[thick] (p9) .. controls (p10) and (p11) .. (p12);
\draw[thick] (p12) .. controls (p13) and (p14) .. (p15);
\draw[thick] (p15) .. controls (p16) and (p17) .. (p18);
\draw[thick] (p18) .. controls (p19) and (p20) .. (p0);
\draw[dashed] (p0) -- (p1) -- (p2) -- (p3);
\draw[dashed] (p3) -- (p4) -- (p5) -- (p6);
\draw[dashed] (p6) -- (p7) -- (p8) -- (p9);
\draw[dashed] (p9) -- (p10) -- (p11) -- (p12);
\draw[dashed] (p12) -- (p13) -- (p14) -- (p15);
\draw[dashed] (p15) -- (p16) -- (p17) -- (p18);
\draw[dashed] (p18) -- (p19) -- (p20) -- (p0);
\end{tikzpicture}
\caption{A closed piecewise B\'ezier curve composed of seven cubic patches.}
\label{piecewiseBeziercurve}
\end{figure}
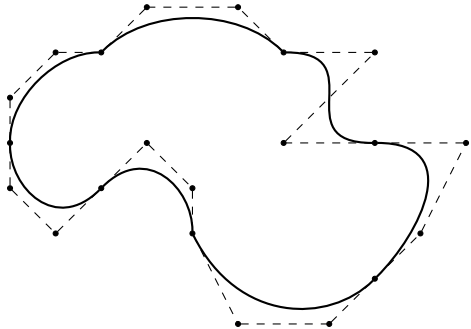
\end{remark}

\begin{remark}
In this paper, since each B\'ezier patch has the same degree $d$, the curve is said to be \textit{uniform in degree}. Nevertheless one can easily build piecewise B\'ezier curves with patches of different degrees. 
\end{remark}


Adapting the proof of the classical Stone-Weierstrass theorem, one can prove the following result (which corresponds to a particular case of the classical Bishop theorem, see \cite{bishop1961}).

\begin{theorem}
Let $f \in \mathcal{C}([0,1],\mathbb{R}^2)$. For all $\varepsilon > 0$ and all $d \in \mathbb{N}^*$, there exist $N \in \mathbb{N}^*$ and a set of $N(d+1)$ control points $P_{1,0},\ldots, P_{1,d},\ldots, P_{N,d}$, satisfying the continuity relations, such that $\| f(t) - B([P_{1,0},\ldots,P_{N,d}],t) \|_{\R^2} \leq \varepsilon$ for all $t \in [0,1]$.
\end{theorem}

This result fully justifies the use of piecewise B\'ezier curves in order to approximate two-dimensional bounded shapes.

\begin{remark}
Recall that the use of polar coordinates, where the radius is expanded in a truncated Fourier series, is another common and efficient strategy in order to approximate two-dimensional shapes (see, e.g., \cite{Dambrine} in the context of inclusions detection). However it has two main drawbacks. Firstly it allows to represent only star-shaped domains and secondly, due to a classical oscillation phenomenon, it cannot represent rigorously straight lines (see, e.g., \cite[Figure~5 p.140]{CauDam12} in the context of inclusions detection). The use of piecewise B\'ezier curves is then an alternative in order to approximate non star-shaped domains and straight lines (see Section~\ref{SectionNumericalSimu} for some numerical simulations in the context of inclusions detection). To conclude this remark, let us recall that the \textit{flip procedure}, which is the main feature of this paper, is based on the detection of potential collisions between two parts of the boundary of the approximated shape (see Section~\ref{SectionFlip} for more details). Thus, it is worth precising that a parameterization based on polar coordinates, where the radius is expanded in a truncated Fourier series, is not adapted to prevent such collisions, in contrary to a piecewise B\'ezier parameterization (see Section~\ref{intersectdetection} for details).
\end{remark}

\section{Intersecting control polygons detection and flip procedure} \label{SectionFlip}
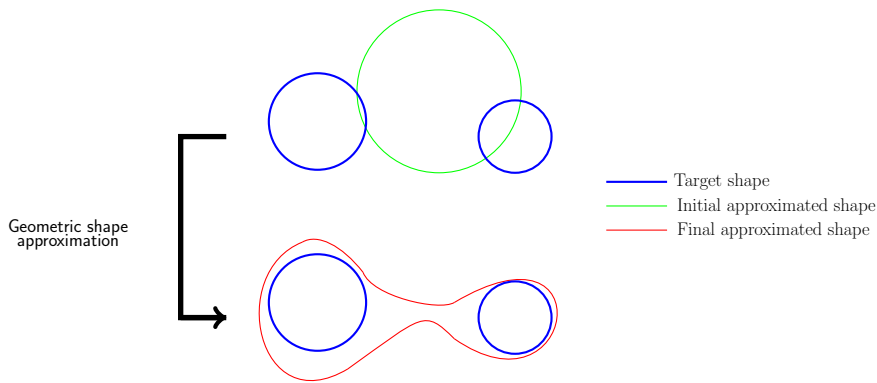
\begin{figure}[h]
\centering
\begin{tikzpicture}[scale=0.4, every node/.style={scale=0.4}]
\draw[green] (1.5,7) circle (2.7cm);
\draw[blue,thick] (4,5.5) circle (1.2cm);
\draw[blue,thick] (-2.5,6) circle (1.6cm);
\draw[line width=2pt] (-5.5,5.5) -- (-7,5.5) -- (-7,-0.5) -- (-5.5,-0.5);
\draw[line width=2pt,->] (-7,-0.5) -- (-5.5,-0.5);
\node (algo) at (-10.7,2.5) {\LARGE \textcolor{black}{\textsf{Geometric shape}}};
\node (algosuite) at (-10.7,2) {\LARGE \textcolor{black}{\textsf{approximation}}};
\draw[red] (2,0) .. controls (1.5,-0.3) and (-0.7,0.2) .. (-1,1);
\draw[red] (-1,1) .. controls (-1.7,1.9) and (-2.5,2.3) .. (-3,2);
\draw[red] (-3,2) .. controls (-5.5,1) and (-4.5,-4) .. (-1.5,-2.2);
\draw[red] (-1.5,-2.2) .. controls (1,-0.4) and (1,-0.2) .. (2,-1.2);
\draw[red] (2,-1.2) .. controls (6,-4) and (7,3) .. (2,0);
\draw[blue,thick] (4,-0.5) circle (1.2cm);
\draw[blue,thick] (-2.5,0) circle (1.6cm);
\draw[blue,thick] (7,4) -- (9,4);
\node (target) at (10.8,4) {\LARGE Target shape};
\draw[green] (7,3.2) -- (9,3.2);
\node (init) at (12.6,3.2) {\LARGE Initial approximated shape};
\draw[red] (7,2.4) -- (9,2.4);
\node (final) at (12.5,2.4) {\LARGE Final approximated shape};
\end{tikzpicture}
\caption{A geometric shape approximation of a two-component target shape starting from a one-component initial approximated shape. The final approximated shape surrounds the two components.}
\label{casFabien}
\end{figure}

In this paper we are interested in geometric two-dimensional shape approximation problems in which the target shape can have multiple connected components but the number of components is unknown. In such a case, starting a classical geometric approximation with a one-component initial shape may lead to the situation depicted in Figure~\ref{casFabien}, that is, the deformation flow makes the boundary evolve until it surrounds all the components of the target shape. This classical phenomenon tends to create a collision between two parts of the boundary of the approximated shape. 

In this paper our major aim is to provide a simple and new concept (called \textit{flip procedure}) that can be added to any shape approximation algorithm based on piecewise B\'ezier curves, 
and which allows to change the topology of the approximated shape. Precisely, the flip procedure allows to divide a one-component shape into a two-component shape.

\begin{remark}
In this paper, we focus on piecewise cubic B\'ezier curve ($d=3$, see Remark~\ref{remarkd3}). However, this method could be easily extended to any $d \geq 2$.
\end{remark}

\subsection{Overview}
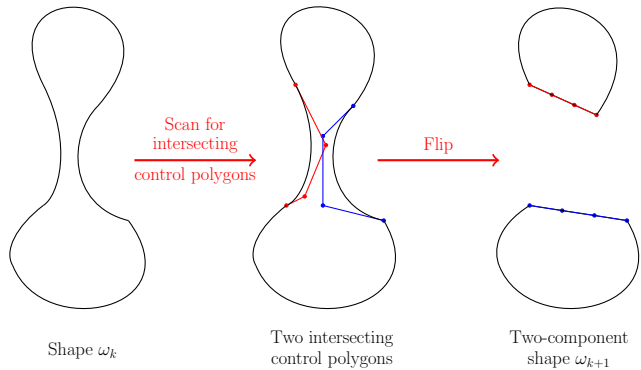
\begin{figure}[h]
\centering
\begin{tikzpicture}[scale=0.4, every node/.style={scale=0.4}]
\draw (-9,3) .. controls (-8,1) and (-8.5,-0.7) .. (-9,-1);
\draw (-9,-1) .. controls (-9.9,-1.7) and (-10.3,-2.5) .. (-10,-3);
\draw (-10,-3) .. controls (-9,-5.5) and (-4,-4.5) .. (-6.2,-1.5);
\draw (-6.2,-1.5) .. controls (-8.2,-1) and (-8.2,1.3) .. (-7.2,2.3);
\draw (-7.2,2.3) .. controls (-4,6) and (-11,7) .. (-9,3);
\node (gamma_k) at (-7.7,-5.8) {\LARGE Shape $\omega_k$};

\node (flip1) at (-4,1.7) {\textcolor{red}{\LARGE Scan for}};
\node (flip2) at (-4,1) {\textcolor{red}{\LARGE intersecting}};
\node (flip3) at (-4,0) {\textcolor{red}{\LARGE control polygons}};
\draw[red,thick,->] (-6,0.5) -- (-2,0.5);

\node[circle,fill=red,draw=red,scale=0.3] (P0) at (-0.7,3) {};
\node[circle,fill=red,draw=red,scale=0.3] (P1) at (0.3,1) {};
\node[circle,fill=red,draw=red,scale=0.3] (P2) at (-0.4,-0.7) {};
\node[circle,fill=red,draw=red,scale=0.3] (P3) at (-1,-1) {};
\node[circle,fill=blue,draw=blue,scale=0.3] (Q0) at (2.2,-1.5) {};
\node[circle,fill=blue,draw=blue,scale=0.3] (Q1) at (0.2,-1) {};
\node[circle,fill=blue,draw=blue,scale=0.3] (Q2) at (0.2,1.3) {};
\node[circle,fill=blue,draw=blue,scale=0.3] (Q3) at (1.2,2.3) {};
\draw[red] (P0) -- (P1) -- (P2) -- (P3);
\draw[blue] (Q0) -- (Q1) -- (Q2) -- (Q3);
\draw (P0) .. controls (P1) and (P2) .. (P3);
\draw (-1,-1) .. controls (-1.9,-1.7) and (-2.3,-2.5) .. (-2,-3);
\draw (-2,-3) .. controls (-1,-5.5) and (4,-4.5) .. (2.2,-1.5);
\draw (Q0) .. controls (Q1) and (Q2) .. (Q3);
\draw (Q3) .. controls (4,6) and (-3,7) .. (P0);
\node (intersection) at (0.5,-5.4) {\LARGE Two intersecting};
\node (intersection2) at (0.5,-6.1) {\LARGE control polygons};

\node (flip1) at (4,1) {\textcolor{red}{\LARGE Flip}};
\draw[red,thick,->] (2,0.5) -- (6,0.5);

\node[circle,fill=red,draw=red,scale=0.3] (PP0) at (7,3) {};
\node[circle,fill=red,draw=red,scale=0.3] (PP1) at (8.47,2.33) {};
\node[circle,fill=blue,draw=blue,scale=0.3] (PP2) at (8.07,-1.17) {};
\node[circle,fill=blue,draw=blue,scale=0.3] (PP3) at (7,-1) {};
\node[circle,fill=blue,draw=blue,scale=0.3] (QQ0) at (10.2,-1.5) {};
\node[circle,fill=blue,draw=blue,scale=0.3] (QQ1) at (9.13,-1.33) {};
\node[circle,fill=red,draw=red,scale=0.3] (QQ2) at (7.73,2.67) {};
\node[circle,fill=red,draw=red,scale=0.3] (QQ3) at (9.2,2) {};

\draw (7,-1) .. controls (5.9,-1.7) and (5.8,-2.5) .. (6,-3);
\draw (6,-3) .. controls (7,-5.5) and (12,-4.5) .. (10.2,-1.5);
\draw (QQ0) .. controls (QQ1) and (PP2) .. (PP3); 
\draw (9.2,2) .. controls (12,6) and (5,7) .. (7,3);
\draw (PP0) .. controls (QQ2) and (PP1) .. (QQ3); 
\draw[red] (PP0) -- (QQ2) -- (PP1) -- (QQ3);
\draw[blue] (QQ0) -- (QQ1) -- (PP2) -- (PP3);
\node (flipped) at (8.3,-5.4) {\LARGE Two-component};
\node (flipped2) at (8.3,-6.1) {\LARGE shape $\omega_{k+1}$};
\end{tikzpicture}
\caption{Overview of the complete procedure.}
\label{detectionandflip}
\end{figure}

Let us consider a general geometric shape approximation algorithm in which the boundary of the approximated shape is parameterized by a piecewise cubic B\'ezier curve. It starts from a one-component initial shape $\omega_0$ and produces a sequence of one-component shapes $(\omega_k)_{k \geq 0}$ by deforming the boundary at each step. Our idea consists in two phases (that are summarized in Figure~\ref{detectionandflip}):
\begin{enumerate}
\item check, at each step of the approximation algorithm, if the current shape $\omega_k$ is in the situation depicted in Figure~\ref{casFabien}, that is, if two parts of the boundary are very close to each other. The parameterization by piecewise B\'ezier curves allows us to prevent such a situation by looking for intersecting control polygons. This procedure will be called \textit{intersecting control polygons detection} and will be detailed in Section~\ref{intersectdetection};
\item if some control polygons intersect each other, we apply the flip procedure in order to obtain a two-component shape by keeping unchanged all other control polygons. The flip procedure is detailed in Section~\ref{flipalgo}.
\end{enumerate}

{ For the sake of simplicity of presentation, we will assume that, during the evolution of the shape~$\omega_k$, situations of intersection between control polygons have always the same pattern:}
\begin{enumerate}
\item[(A1)] { only two situations of intersection between control polygons are possible: either one control polygon intersects exactly another one, or one control polygon intersects exactly two consecutive ones (see Figure~\ref{2casintersections});}

\begin{figure}[h]
\centering
\subfigure[Case with two control polygons]{
\begin{tikzpicture}[scale=0.8, every node/.style={scale=0.7}]
\node[circle,fill=red,draw=red,scale=0.3] (P0) at (-0.7,3) {};
\node[red] at (-1.2,3) {$P_{0}$};
\node[circle,fill=red,draw=red,scale=0.3] (P1) at (0.3,1) {};
\node[red] at (1,1) {$P_{1}$};
\node[circle,fill=red,draw=red,scale=0.3] (P2) at (-0.4,-0.7) {};
\node[red] at (-0.5,-1.2) {$P_{2}$};
\node[circle,fill=red,draw=red,scale=0.3] (P3) at (-1,-1) {};
\node[red] at (-1.5,-1) {$P_{3}$};
\node[circle,fill=blue,draw=blue,scale=0.3] (Q0) at (2.2,-1.5) {};
\node[blue] at (2.7,-1.5) {$Q_{0}$};
\node[circle,fill=blue,draw=blue,scale=0.3] (Q1) at (0.2,-1) {};
\node[blue] at (0.2,-1.4) {$Q_{1}$};
\node[circle,fill=blue,draw=blue,scale=0.3] (Q2) at (0.2,1.3) {};
\node[blue] at (-0.5,1.3) {$Q_{2}$};
\node[circle,fill=blue,draw=blue,scale=0.3] (Q3) at (1.2,2.3) {};
\node[blue] at (1.7,2.3) {$Q_{3}$};
\draw[red,thick] (P0) -- (P1) -- (P2) -- (P3);
\draw[blue,thick] (Q0) -- (Q1) -- (Q2) -- (Q3);
\draw[red] (P0) .. controls (P1) and (P2) .. (P3);
\draw[blue] (Q0) .. controls (Q1) and (Q2) .. (Q3);
\end{tikzpicture}
\label{ctrlpolygonsintersection}
}
\hspace{2cm}
\subfigure[Case with three control polygons]{
\begin{tikzpicture}[scale=0.45, every node/.style={scale=0.4}]
\node[circle,fill=red,draw=red,scale=0.5] (P0) at (0,8) {};
\node[circle,fill=red,draw=red,scale=0.5] (P1) at (3,5.5) {};
\node[circle,fill=red,draw=red,scale=0.5] (P2) at (3.3,2) {};
\node[circle,fill=red,draw=red,scale=0.5] (P3) at (0,0) {};
\node[red] (P_0) at (0.2,8.4) {$P_{0}$};
\node[red] (P_1) at (3.3,5.8) {$P_{1}$};
\node[red] (P_2) at (3.6,1.7) {$P_{2}$};
\node[red] (P_3) at (0.3,-0.3) {$P_{3}$};
\draw[red,thick] (P0) -- (P1) -- (P2) -- (P3);
\draw[red] (P0) .. controls (P1) and (P2) .. (P3);

\node[circle,fill=blue,draw=blue,scale=0.5] (Q0) at (8,0) {};
\node[circle,fill=blue,draw=blue,scale=0.5] (Q1) at (4.5,0.5) {};
\node[circle,fill=blue,draw=blue,scale=0.5] (Q2) at (3,1.2) {};
\node[circle,fill=blue,draw=blue,scale=0.5] (Q3) at (2.7,3) {};
\node[blue] (Q_0) at (8.4,0) {$Q_{0}$};
\node[blue] (Q_1) at (4.2,0.2) {$Q_{1}$};
\node[blue] (Q_2) at (3,0.7) {$Q_{2}$};
\node[blue] (Q_3) at (3.3,3) {$Q_{3}$};
\draw[blue,thick] (Q0) -- (Q1) -- (Q2) -- (Q3);
\draw[blue] (Q0) .. controls (Q1) and (Q2) .. (Q3);

\node[circle,fill=green,draw=blue,scale=0.5] (R0) at (2.7,3) {};
\node[circle,fill=green,draw=green,scale=0.5] (R1) at (3,6.5) {};
\node[circle,fill=green,draw=green,scale=0.5] (R2) at (5,7.5) {};
\node[circle,fill=green,draw=green,scale=0.5] (R3) at (7,6) {};
\node[green] (R_0) at (2.2,3) {$R_{0}$};
\node[green] (R_1) at (3,6.8) {$R_{1}$};
\node[green] (R_2) at (5,7.8) {$R_{2}$};
\node[green] (R_3) at (7.3,6.3) {$R_{3}$};
\draw[green,thick] (R0) -- (R1) -- (R2) -- (R3);
\draw[green] (R0) .. controls (R1) and (R2) .. (R3);
\end{tikzpicture}
\label{cas3patches}
}
\caption{Assumption~(A1).}
\label{2casintersections}
\end{figure}
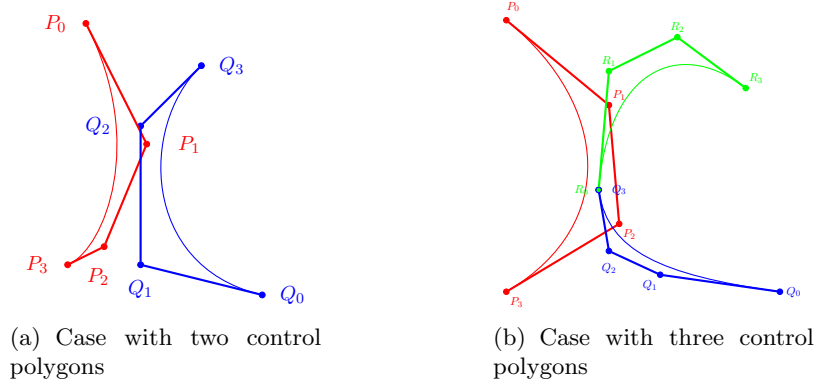

\item[(A2)] { futhermore, at each iteration, at most one situation of intersection occurs.}
\end{enumerate}
{ Assumptions (A1)-(A2) are ordinarily satisfied in practice, in particular in all numerical simulations we made (see Section~\ref{SectionNumericalSimu}). Removing these assumptions should not involve new deep ideas or new concept, however the complete algorithmic description and implementation would become considerably tricky and challenging. It is not our aim to deal with this issue in this paper.}

\begin{remark}
{ In Figure~\ref{2casintersections}, note that collisions between two patches are not excluded. In that case, we say that the shape $\omega_k$ is \textit{self-intersecting}. However, from Assumption (A3) enunciated later, this situation does not jeopardize the integrity of Algorithm~$\mathcal{A}$ presented in Section~\ref{secalgo} (see in particular Step~\ref{flipoupasflip} of the algorithm).}
\end{remark}

\begin{remark}[Controlling the size of the patches using split and merge functions]\label{remcontrolsize}
In order to maintain numerical stability, one should control the size of the control polygons (that is, the diameter of their convex hull) in a range $[S_{\text{min}},S_{\text{max}}]$ with $0 < S_{\text{min}} < S_{\text{max}}$. This avoids to deal with very large patches and/or very small ones.
To this end, the diameter of the convex hull of each control polygon can be computed at each iteration. If a control polygon does not satisfy the size condition, it is either split into two control polygons or merged with a neighbor one. The split and merge functions (see \cite{Ruatta:Freeform} for more details) are inverse operations and both use interpolation in order to compute the new control polygons (see Figure~\ref{splitandmerge}). The split function divides a control polygon into two. Precisely, it interpolates the first half of the patch and, in a second time, interpolates the other half. Since each half of the patch is a B\'ezier curve, the shape is not modified after a split. The merge function is the reverse operation. From two consecutive control polygons $\mathbf{Q}$ and $\mathbf{R}$, it computes one patch that interpolates the four points $B(\mathbf{Q},0)$, $B(\mathbf{Q},\frac{2}{3})$, $B(\mathbf{R},\frac{1}{3})$ and $B(\mathbf{R},1)$. Then one starts from seven control points and ends with four. Note that merging polygons modifies slightly the boundary.
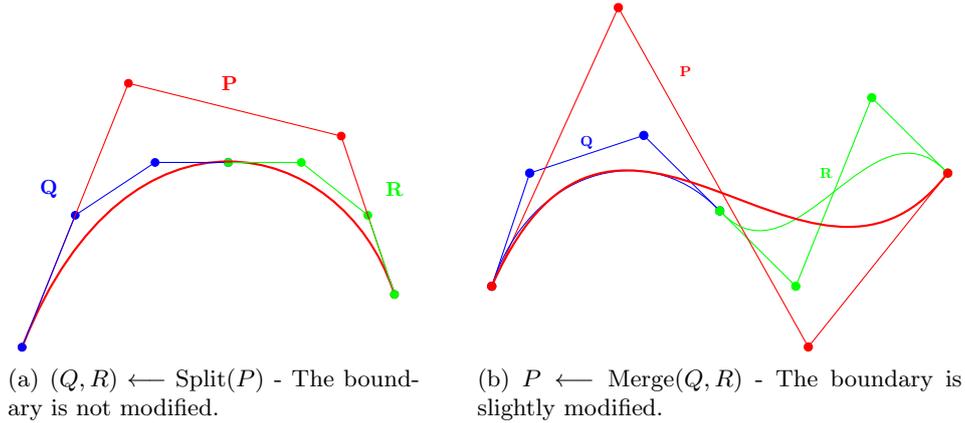
\begin{figure}[h]
\centering
\subfigure[$(Q,R) \longleftarrow \text{Split}(P)$ - The boundary is not modified.]{
\begin{tikzpicture}[scale=0.7, every node/.style={scale=0.7}]
\node[circle,fill=red,draw=red,scale=0.4] (p0) at (0,0) {};
\node[circle,fill=red,draw=red,scale=0.4] (p1) at (2,5) {};
\node[circle,fill=red,draw=red,scale=0.4] (p2) at (6,4) {};
\node[circle,fill=red,draw=red,scale=0.4] (p3) at (7,1) {};
\node[circle,fill=blue,draw=blue,scale=0.4] (q0) at (0,0) {};
\node[circle,fill=blue,draw=blue,scale=0.4] (q1) at (1,2.5) {};
\node[circle,fill=blue,draw=blue,scale=0.4] (q2) at (2.5,3.5) {};
\node[circle,fill=blue,draw=blue,scale=0.4] (q3) at (3.875,3.5) {};
\node[circle,fill=green,draw=green,scale=0.4] (r0) at (3.875,3.5) {};
\node[circle,fill=green,draw=green,scale=0.4] (r1) at (5.25,3.5) {};
\node[circle,fill=green,draw=green,scale=0.4] (r2) at (6.5,2.5) {};
\node[circle,fill=green,draw=green,scale=0.4] (r3) at (7,1) {};
\draw[red,thick] (p0) .. controls (p1) and (p2) .. (p3);
\draw[red] (p0) -- (p1) -- (p2) -- (p3);
\draw[blue] (q0) -- (q1) -- (q2) -- (q3);
\draw[green] (r0) -- (r1) -- (r2) -- (r3);
\node[red] (P) at (3.9,5) {$\mathbf{P}$};
\node[blue] (Q) at (0.5,3) {$\mathbf{Q}$};
\node[green] (R) at (7,3) {$\mathbf{R}$};
\end{tikzpicture}
\label{split}
}
\hspace{5mm}
\subfigure[$P \longleftarrow \text{Merge}(Q,R)$ - The boundary is slightly modified.]{
\begin{tikzpicture}[scale=0.5, every node/.style={scale=0.5}]
\node[circle,fill=blue,draw=blue,scale=0.6] (Q0) at (0,0) {};
\node[circle,fill=blue,draw=blue,scale=0.6] (Q1) at (1,3) {};
\node[circle,fill=blue,draw=blue,scale=0.6] (Q2) at (4,4) {};
\node[circle,fill=blue,draw=blue,scale=0.6] (Q3) at (6,2) {};
\node[circle,fill=green,draw=green,scale=0.6] (R0) at (6,2) {};
\node[circle,fill=green,draw=green,scale=0.6] (R1) at (8,0) {};
\node[circle,fill=green,draw=green,scale=0.6] (R2) at (10,5) {};
\node[circle,fill=green,draw=green,scale=0.6] (R3) at (12,3) {};
\node[circle,fill=red,draw=red,scale=0.6] (P0) at (0,0) {};
\node[circle,fill=red,draw=red,scale=0.6] (P1) at (3.33,7.39) {};
\node[circle,fill=red,draw=red,scale=0.6] (P2) at (8.33,-1.61) {};
\node[circle,fill=red,draw=red,scale=0.6] (P3) at (12,3) {};
\draw[blue] (Q0) -- (Q1) -- (Q2) -- (Q3);
\draw[green] (R0) -- (R1) -- (R2) -- (R3);
\draw[red] (P0) -- (P1) -- (P2) -- (P3);
\draw[blue] (Q0) .. controls (Q1) and (Q2) .. (Q3);
\draw[green] (R0) .. controls (R1) and (R2) .. (R3);
\draw[red,thick] (P0) .. controls (P1) and (P2) .. (P3);
\node[blue] (Q) at (2.5,3.8) {$\mathbf{Q}$};
\node[green] (R) at (8.8,3) {$\mathbf{R}$};
\node[red] (P) at (5.1,5.7) {$\mathbf{P}$};
\end{tikzpicture}
\label{merge}
}
\caption{Examples of the split and merge functions.}
\label{splitandmerge}
\end{figure}
\end{remark}

\subsection{Intersecting control polygons detection}\label{intersectdetection}
\begin{figure}[h]
\centering
\begin{tikzpicture}[scale=0.6, every node/.style={scale=0.9}]

\draw[->] (-0.5,0.5) -- (0.5,0.5);
\draw[->] (-0.5,0.5) -- (-0.5,1.5);
\node (x) at (0,0.2) {$x$};
\node (y) at (-0.8,1) {$y$};
\draw[red,thick,fill=red!20] (0,2) -- (2,2) -- (2,4) -- (0,4) -- (0,2);
\draw[red,thick,fill=red!20] (2,2) -- (4,2) -- (4,4) -- (2,4) -- (2,2);
\draw[red,thick,fill=red!20] (4,0) -- (8,0) -- (8,2) -- (4,2) -- (4,0);
\draw[red,thick,fill=red!20] (8,1) -- (10,1) -- (10,4) -- (8,4) -- (8,1);
\draw[red,thick,fill=red!20] (6,4) -- (8,4) -- (8,6) -- (6,6) -- (6,4);
\draw[red,thick,fill=red!20] (2,6) -- (6,6) -- (6,7) -- (2,7) -- (2,6);
\draw[red,thick,fill=red!20] (0,4) -- (2,4) -- (2,6) -- (0,6) -- (0,4);
\node[circle,fill=black,draw=black,scale=0.3] (p0) at (0,4) {};
\node[circle,fill=black,draw=black,scale=0.3] (p1) at (0,3) {};
\node[circle,fill=black,draw=black,scale=0.3] (p2) at (1,2) {};
\node[circle,fill=black,draw=black,scale=0.3] (p3) at (2,3) {};
\node[circle,fill=black,draw=black,scale=0.3] (p4) at (3,4) {};
\node[circle,fill=black,draw=black,scale=0.3] (p5) at (4,3) {};
\node[circle,fill=black,draw=black,scale=0.3] (p6) at (4,2) {};
\node[circle,fill=black,draw=black,scale=0.3] (p7) at (5,0) {};
\node[circle,fill=black,draw=black,scale=0.3] (p8) at (7,0) {};
\node[circle,fill=black,draw=black,scale=0.3] (p9) at (8,1) {};
\node[circle,fill=black,draw=black,scale=0.3] (p10) at (9,2) {};
\node[circle,fill=black,draw=black,scale=0.3] (p11) at (10,4) {};
\node[circle,fill=black,draw=black,scale=0.3] (p12) at (8,4) {};
\node[circle,fill=black,draw=black,scale=0.3] (p13) at (6,4) {};
\node[circle,fill=black,draw=black,scale=0.3] (p14) at (8,6) {};
\node[circle,fill=black,draw=black,scale=0.3] (p15) at (6,6) {};
\node[circle,fill=black,draw=black,scale=0.3] (p16) at (5,7) {};
\node[circle,fill=black,draw=black,scale=0.3] (p17) at (3,7) {};
\node[circle,fill=black,draw=black,scale=0.3] (p18) at (2,6) {};
\node[circle,fill=black,draw=black,scale=0.3] (p19) at (1,6) {};
\node[circle,fill=black,draw=black,scale=0.3] (p20) at (0,5) {};
\draw[dashed] (p0) .. controls (p1) and (p2) .. (p3);
\draw[dashed] (p3) .. controls (p4) and (p5) .. (p6);
\draw[dashed] (p6) .. controls (p7) and (p8) .. (p9);
\draw[dashed] (p9) .. controls (p10) and (p11) .. (p12);
\draw[dashed] (p12) .. controls (p13) and (p14) .. (p15);
\draw[dashed] (p15) .. controls (p16) and (p17) .. (p18);
\draw[dashed] (p18) .. controls (p19) and (p20) .. (p0);
\draw[] (p0) -- (p1) -- (p2) -- (p3);
\draw[] (p3) -- (p4) -- (p5) -- (p6);
\draw[] (p6) -- (p7) -- (p8) -- (p9);
\draw[] (p9) -- (p10) -- (p11) -- (p12);
\draw[] (p12) -- (p13) -- (p14) -- (p15);
\draw[] (p15) -- (p16) -- (p17) -- (p18);
\draw[] (p18) -- (p19) -- (p20) -- (p0);
\end{tikzpicture}
\caption{AABBs of control polygons.}
\label{AABB}
\end{figure}
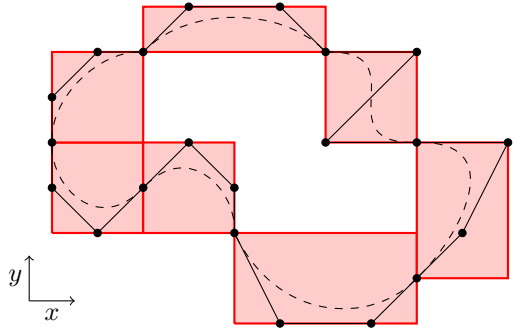

Checking if each control polygon intersects another one may be very expensive in terms of computations. Axis-Aligned Bounding Boxes (AABBs) are a very common tool in Computer Graphics and Computational Geometry in order to detect the collision of two objects (see, e.g., \cite{Ericson}), with a relatively low computational cost. AABB is defined as the smallest rectangle, whose sides are aligned with the axes, containing the control polygon (see Figure~\ref{AABB}). 

A necessary condition for two intersecting control polygons is clearly the intersection of their respective AABBs. As a consequence, instead of looking directly for intersecting control polygons, we first look for intersecting AABBs. Thus, the intersecting control polygons detection consists in two steps:
\begin{enumerate}
\item we first list all the pairs of intersecting AABBs;
\item in a second time, we check these pairs in order to see if the associated control polygons intersect. To do so, we directly check the $9$ segment-segment intersections of the polygons (see, e.g., \cite[p.~28-30]{O'Rourke:1998:CGC:521378}).
\end{enumerate}
Finally, each pair of intersecting control polygons will be given as input to the flip procedure detailed in the following section.

\subsection{The flip procedure}\label{flipalgo}

From Assumption~(A1), only two cases of intersecting control polygons are considered (see Figure~\ref{2casintersections}). The flip procedure described in this section is a simple tool that can be easily implemented and that handles these two situations.

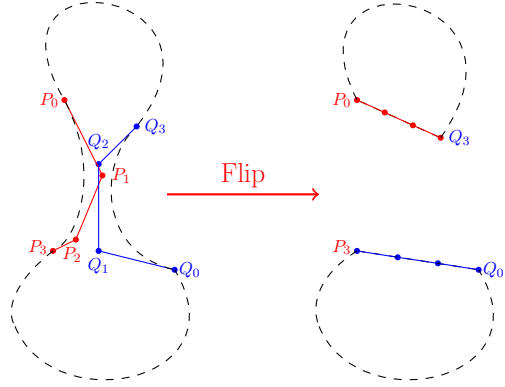
\begin{figure}[h]
\begin{tikzpicture}[scale=0.5, every node/.style={scale=0.6}]
\node[circle,fill=red,draw=red,scale=0.3] (P0) at (-0.7,3) {};
\node[circle,fill=red,draw=red,scale=0.3] (P1) at (0.3,1) {};
\node[circle,fill=red,draw=red,scale=0.3] (P2) at (-0.4,-0.7) {};
\node[circle,fill=red,draw=red,scale=0.3] (P3) at (-1,-1) {};
\node[red] (p0) at (-1.1,3) {$P_0$};
\node[red] (p1) at (0.8,1) {$P_1$};
\node[red] (p2) at (-0.5,-1.1) {$P_2$};
\node[red] (p3) at (-1.4,-0.9) {$P_3$};
\node[circle,fill=blue,draw=blue,scale=0.3] (Q0) at (2.2,-1.5) {};
\node[circle,fill=blue,draw=blue,scale=0.3] (Q1) at (0.2,-1) {};
\node[circle,fill=blue,draw=blue,scale=0.3] (Q2) at (0.2,1.3) {};
\node[circle,fill=blue,draw=blue,scale=0.3] (Q3) at (1.2,2.3) {};
\node[blue] (q0) at (2.6,-1.5) {$Q_0$};
\node[blue] (q1) at (0.2,-1.4) {$Q_1$};
\node[blue] (q2) at (0.2,1.9) {$Q_2$};
\node[blue] (q3) at (1.7,2.3) {$Q_3$};

\draw[red] (P0) -- (P1) -- (P2) -- (P3);
\draw[blue] (Q0) -- (Q1) -- (Q2) -- (Q3);
\draw[dashed] (P0) .. controls (P1) and (P2) .. (P3);
\draw[dashed] (-1,-1) .. controls (-1.9,-1.7) and (-2.3,-2.5) .. (-2,-3);
\draw[dashed] (-2,-3) .. controls (-1,-5.5) and (4,-4.5) .. (2.2,-1.5);
\draw[dashed] (Q0) .. controls (Q1) and (Q2) .. (Q3);
\draw[dashed] (Q3) .. controls (4,6) and (-3,7) .. (P0);

\node (flip1) at (4,1) {\textcolor{red}{\LARGE Flip}};
\draw[red,thick,->] (2,0.5) -- (6,0.5);

\node[circle,fill=red,draw=red,scale=0.3] (PP0) at (7,3) {};
\node[circle,fill=red,draw=red,scale=0.3] (PP1) at (8.47,2.33) {};
\node[circle,fill=blue,draw=blue,scale=0.3] (PP2) at (8.07,-1.17) {};
\node[circle,fill=blue,draw=blue,scale=0.3] (PP3) at (7,-1) {};
\node[circle,fill=blue,draw=blue,scale=0.3] (QQ0) at (10.2,-1.5) {};
\node[circle,fill=blue,draw=blue,scale=0.3] (QQ1) at (9.13,-1.33) {};
\node[circle,fill=red,draw=red,scale=0.3] (QQ2) at (7.73,2.67) {};
\node[circle,fill=red,draw=red,scale=0.3] (QQ3) at (9.2,2) {};

\node[red] (pp0) at (6.6,3) {$P_0$};
\node[red] (pp3) at (6.6,-0.9) {$P_3$};

\node[blue] (qq0) at (10.6,-1.5) {$Q_0$};
\node[blue] (qq3) at (9.7,2) {$Q_3$};

\draw[dashed] (7,-1) .. controls (5.9,-1.7) and (5.8,-2.5) .. (6,-3);
\draw[dashed] (6,-3) .. controls (7,-5.5) and (12,-4.5) .. (10.2,-1.5);
\draw[dashed] (QQ0) .. controls (QQ1) and (PP2) .. (PP3); 
\draw[dashed] (9.2,2) .. controls (12,6) and (5,7) .. (7,3);
\draw[dashed] (PP0) .. controls (QQ2) and (PP1) .. (QQ3); 
\draw[red] (PP0) -- (QQ2) -- (PP1) -- (QQ3);
\draw[blue] (QQ0) -- (QQ1) -- (PP2) -- (PP3);

\end{tikzpicture}
\caption{Flip procedure - Case of two control polygons.}
\label{flipCas1}
\end{figure}

\paragraph{First case: two intersecting polygons.}
From $\mathbf{P}=\{P_0,P_1,P_2,P_3\}$ and $\mathbf{Q}=\{Q_0,Q_1,Q_2,Q_3\}$ being two intersecting polygons of a same connected component, the flip procedure builds two new polygons as follows, (see Figure~\ref{flipCas1}):
$$ \left\lbrace P_0,P_0+\frac{1}{3}\overrightarrow{P_0Q_3},P_0+\frac{2}{3}\overrightarrow{P_0Q_3},Q_3 \right\rbrace \qquad \text{and} \qquad \left\lbrace Q_0,Q_0+\frac{1}{3}\overrightarrow{Q_0P_3},Q_0+\frac{2}{3}\overrightarrow{Q_0P_3},P_3 \right\rbrace. $$

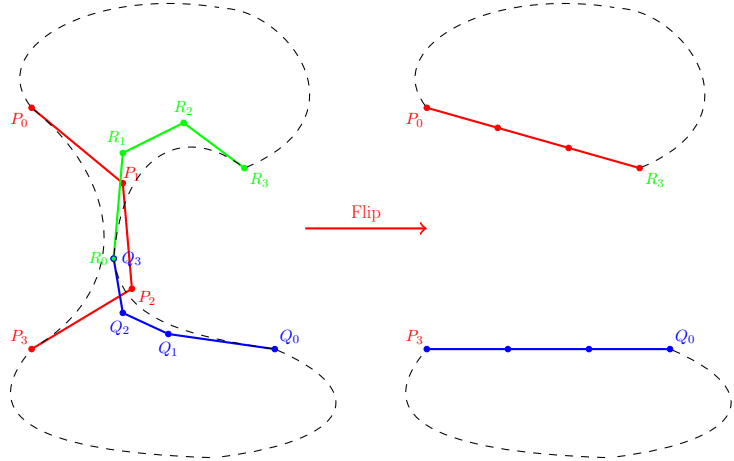
\begin{figure}[h]
\centering
\begin{tikzpicture}[scale=0.4, every node/.style={scale=0.4}]
\node[circle,fill=red,draw=red,scale=0.5] (P0) at (0,8) {};
\node[circle,fill=red,draw=red,scale=0.5] (P1) at (3,5.5) {};
\node[circle,fill=red,draw=red,scale=0.5] (P2) at (3.3,2) {};
\node[circle,fill=red,draw=red,scale=0.5] (P3) at (0,0) {};
\node[red,scale = 1.4] (P_0) at (-0.4,7.6) {$P_{0}$};
\node[red,scale = 1.4] (P_1) at (3.3,5.8) {$P_{1}$};
\node[red,scale = 1.4] (P_2) at (3.8,1.7) {$P_{2}$};
\node[red,scale = 1.4] (P_3) at (-0.4,0.4) {$P_{3}$};
\draw[red,thick] (P0) -- (P1) -- (P2) -- (P3);
\draw[dashed] (P0) .. controls (P1) and (P2) .. (P3);

\node[circle,fill=blue,draw=blue,scale=0.5] (Q0) at (8,0) {};
\node[circle,fill=blue,draw=blue,scale=0.5] (Q1) at (4.5,0.5) {};
\node[circle,fill=blue,draw=blue,scale=0.5] (Q2) at (3,1.2) {};
\node[circle,fill=blue,draw=blue,scale=0.5] (Q3) at (2.7,3) {};
\node[blue,scale = 1.4] (Q_0) at (8.5,0.4) {$Q_{0}$};
\node[blue,scale = 1.4] (Q_1) at (4.5,0) {$Q_{1}$};
\node[blue,scale = 1.4] (Q_2) at (2.9,0.7) {$Q_{2}$};
\node[blue,scale = 1.4] (Q_3) at (3.3,3) {$Q_{3}$};
\draw[blue,thick] (Q0) -- (Q1) -- (Q2) -- (Q3);
\draw[dashed] (Q0) .. controls (Q1) and (Q2) .. (Q3);

\node[circle,fill=green,draw=blue,scale=0.5] (R0) at (2.7,3) {};
\node[circle,fill=green,draw=green,scale=0.5] (R1) at (3,6.5) {};
\node[circle,fill=green,draw=green,scale=0.5] (R2) at (5,7.5) {};
\node[circle,fill=green,draw=green,scale=0.5] (R3) at (7,6) {};
\node[green,scale = 1.4] (R_0) at (2.2,3) {$R_{0}$};
\node[green,scale = 1.4] (R_1) at (2.8,7) {$R_{1}$};
\node[green,scale = 1.4] (R_2) at (5,8) {$R_{2}$};
\node[green,scale = 1.4] (R_3) at (7.5,5.6) {$R_{3}$};
\draw[green,thick] (R0) -- (R1) -- (R2) -- (R3);
\draw[dashed] (R0) .. controls (R1) and (R2) .. (R3);

\node (F1) at (-1.5,10) {};
\node (F2) at (2,12) {};
\node (F3) at (6.3,11.3) {};
\node (G1) at (8.5,11.1) {};
\node (G2) at (11,7.7) {};
\draw[dashed] (P0) .. controls (F1) and (F2) .. (F3);
\draw[dashed] (F3) .. controls (G1) and (G2) .. (R3);
\node (S1) at (11,-1.3) {};
\node (S2) at (10.8,-3) {};
\node (S3) at (6,-3.6) {};
\node (T1) at (0,-3.4) {};
\node (T2) at (-1.8,-2.2) {};
\draw[dashed] (Q0) .. controls (S1) and (S2) .. (S3);
\draw[dashed] (S3) .. controls (T1) and (T2) .. (P3);

\node (flip) at (11,4.5) {\textcolor{red}{\LARGE Flip}};
\draw[red,thick,->] (9,4) -- (13,4);

\node[circle,fill=red,draw=red,scale=0.5] (PR0) at (13+0,8) {};
\node[red,scale = 1.4] (PP_0) at (13-0.4,7.6) {$P_{0}$};
\node[circle,fill=red,draw=red,scale=0.5] (PR1) at (13+0+2.333,8-0.666) {};
\node[circle,fill=red,draw=red,scale=0.5] (PR2) at (13+0+4.666,8-1.332) {};
\node[circle,fill=red,draw=red,scale=0.5] (PR3) at (13+7,6) {};
\node[green,scale = 1.4] (RR_3) at (13+7.5,5.6) {$R_{3}$};
\draw[red,thick] (PR0) -- (PR1) -- (PR2) -- (PR3);

\node[circle,fill=blue,draw=blue,scale=0.5] (PQ0) at (13+8,0) {}; 
\node[blue,scale = 1.4] (Q_0) at (13+8.5,0.4) {$Q_{0}$};
\node[circle,fill=blue,draw=blue,scale=0.5] (PQ1) at (13+8-2.666,0) {}; 
\node[circle,fill=blue,draw=blue,scale=0.5] (PQ2) at (13+8-5.332,0) {}; 
\node[circle,fill=blue,draw=blue,scale=0.5] (PQ3) at (13+0,0) {}; 
\node[red,scale = 1.4] (PP_3) at (13-0.4,0.4) {$P_{3}$};
\draw[blue,thick] (PQ0) -- (PQ1) -- (PQ2) -- (PQ3);

\node (F1) at (13-1.5,10) {};
\node (F2) at (13+2,12) {};
\node (F3) at (13+6.3,11.3) {};
\node (G1) at (13+8.5,11.1) {};
\node (G2) at (13+11,7.7) {};
\draw[dashed] (PR0) .. controls (F1) and (F2) .. (F3);
\draw[dashed] (F3) .. controls (G1) and (G2) .. (PR3);
\node (S1) at (13+11,-1.3) {};
\node (S2) at (13+10.8,-3) {};
\node (S3) at (13+6,-3.6) {};
\node (T1) at (13+0,-3.4) {};
\node (T2) at (13-1.8,-2.2) {};
\draw[dashed] (PQ0) .. controls (S1) and (S2) .. (S3);
\draw[dashed] (S3) .. controls (T1) and (T2) .. (PQ3);
\end{tikzpicture}

\caption{Flip procedure - Case of three control polygons.}
\label{flipCas2}
\end{figure}

\paragraph{Second case: three intersecting polygons.} 
The case with three control polygons is very similar. From $\mathbf{P}=\{P_0,P_1,P_2,P_3\}$ being a control polygon intersecting two consecutive ones $\mathbf{Q}=\{Q_0,Q_1,Q_2,Q_3\}$ and $\mathbf{R}=\{R_0,R_1,R_2,R_3\}$, the flip procedure builds two new polygons as follows, (see Figure~\ref{flipCas2}):
$$ \left\lbrace P_0,P_0+\frac{1}{3}\overrightarrow{P_0R_3},P_0+\frac{2}{3}\overrightarrow{P_0R_3},R_3 \right\rbrace \qquad \text{and} \qquad \left\lbrace Q_0,Q_0+\frac{1}{3}\overrightarrow{Q_0P_3},Q_0+\frac{2}{3}\overrightarrow{Q_0P_3},P_3 \right\rbrace. $$

\bigskip

{ For the sake of simplicity of presentation, we will assume that the following hypothesis is satisfied:}
\begin{enumerate}
\item[(A3)] { the flip procedure does not produce intersecting control polygons.}
\end{enumerate}
{ In the same spirit of Assumptions (A1)-(A2), Assumption (A3) is ordinarily satisfied in practice, in particular in all numerical simulations we made (see Section~\ref{SectionNumericalSimu}).}

\section{Application to multiple-inclusion detection} \label{SectionObjectsDetection}

This section focuses on the problem of reconstructing numerically an obstacle $\omega_{\rm ex}$ living in a larger bounded domain $\Omega$ of $\mathbb{R}^2$ from boundary measurements. Our aim is in particular to test the flip procedure introduced in this paper in the case where $\omega_{\rm ex}$ is a two-component obstacle (see Section~\ref{sectionseveralobstacles}). 

In order to solve numerically the above inverse obstacle problem, we will actually consider a shape optimization problem, by minimizing a shape cost functional. In this paper we use the classical geometrical shape optimization approach, based on shape derivatives and on a shape gradient descent method. We refer to the classical books of Henrot \textit{et al.}~\cite{HP} and of Soko{\l}owski \textit{et al.}~\cite{S-Z} for more details on the techniques of shape differentiability.

Let us fix some notations that will be used in this section.  We denote by $\mathrm{L}^p$, $ \mathrm{W}^{m,p} $ and $\mathrm{H}^s$ the usual Lebesgue and Sobolev spaces. We note in bold the vectorial functions and spaces, such as 
$ \boldsymbol{\mathrm{W}}^{m,p} $. Let~$ \Omega $ be a nonempty bounded and connected open set of $\R^2 $ with a $C^{2,1}$ boundary and let $g \in {\mathrm{H}}^{5/2}(\partial\Omega) $ such that $g \neq 0$. We denote by $\boldsymbol{\mathrm{n}}$ the external unit normal to $\partial \Omega$, and for a smooth enough function~$u$, we denote by~$ \partial_{\boldsymbol{\mathrm{n}}} u $ the normal derivative of~$u$. 

Let $ d_0 > 0 $ be fixed (small). In the sequel $ \mathcal{O}_{d_0} $ stands for the set of all open subsets~$ \omega $ strictly included in~$\Omega$, with a~$ C^{2,1}$ boundary, such that the distance $ \mathrm{d}(x,\partial \Omega)$ from~$x$ to the compact $\partial\Omega$ is strictly greater than $d_0$ for all $ x \in \omega $, and such that $\priv{\Omega}{\omega}$ is connected. Finally we also introduce $\Omega_{d_0}$ an open set with a $C^{\infty}$ boundary such that
$$\displaystyle  \left\lbrace x \in \Omega \, ; \; \mathrm{d}(x,\partial \Omega) > d_0 / 2 \right\rbrace \subset   \Omega_{d_0} \subset  \left\lbrace x \in \Omega \, ; \; \mathrm{d}(x,\partial \Omega) > d_0 / 3 \right\rbrace.$$

\subsection{Problem setting}
We focus on the following inverse problem. Assume that an unknown obstacle $\omega_{\rm ex} \in \mathcal{O}_{d_0}$ is located inside $\Omega$. 
We consider hereafter the Laplace equation in $ \Omega \backslash \overline{\omega_{\rm ex}}$ with homogeneous Dirichlet boundary condition on $\partial \omega_{\rm ex}$ and non-homogeneous Dirichlet boundary condition on $\partial \Omega$. Precisely we denote by $u_{\rm ex} \in \HH^1(\priv{\Omega}{\omega_{\rm ex}})$ the unique solution of the problem
\begin{equation} \label{PbExact}
\left\lbrace 
\begin{array}{rcll}
\displaystyle - \Delta u_{\rm ex}  & = & 0  & \mbox{in } \, \Omega \backslash \overline{\omega_{\rm ex}} , \\
\displaystyle u_{ \rm ex} &= & g  & \mbox{on } \, \partial \Omega  , \\
u_{ \rm ex} &= & 0   & \mbox{on } \, \partial \omega_{\rm ex}.
\end{array}
\right.
\end{equation}
Since $g \in {\mathrm{H}}^{5/2}(\partial\Omega) $ and $ \Omega \backslash \overline{\omega_{\rm ex}}$ has a $C^{2,1}$ boundary, note that $u_{\rm ex}$ belongs to $\HH^3(\priv{\Omega}{\omega_{\rm ex}})$. Our main purpose is to reconstruct the unknown shape~$\omega_{\rm ex}$, assuming that a measurement is done on the exterior boundary $\partial\Omega$. Precisely we assume in this paper that we know exactly the value of the measure $ f_{b} \assign \partial_{\boldsymbol{\mathrm{n}}} u_{\rm ex} \in \HH^{3/2} (\partial\Omega)$ on $\partial\Omega$. Thus, for a given nontrivial Cauchy pair $(g,f_{b}) \in \HH^{5/2}(\partial\Omega) \times \HH^{3/2}(\partial\Omega)$, we are interested in the following geometric inverse problem:
\begin{equation}\label{problemInverse}
\begin{array}{c}
\mbox{\textit{find }}  \omega \in \mathcal{O}_{d_0}  \mbox{\textit{ and }}   u \in \HH^1(\priv{\Omega}{\omega})  \cap C^0(\priv{\Omega}{\omega}) \mbox{\textit{ which satisfies the overdetermined system }} \\[5pt]
\left\lbrace 
\begin{array}{rcll}
\displaystyle - \Delta u  & = & 0  & \mbox{in } \, \Omega \backslash \overline{\omega} , \\
\displaystyle u &= & g  & \mbox{on } \, \partial \Omega  , \\
\displaystyle \partial_{\Gn} u &= & f_{b}  & \mbox{on } \, \partial \Omega  , \\
u &= & 0   & \mbox{on } \, \partial \omega .
\end{array}
\right.
\end{array}
\end{equation}
The existence of a solution is trivial since we assume that the measurement $f_b$ is exact. From the classical Holmgren's theorem (see, e.g., \cite{isakov2006inverse}) one can obtain an identifiability result for this inverse problem which claims that the solution is unique. This fundamental question about uniqueness of a solution to the overdetermined problem~\eqref{problemInverse} was deeply studied, see for example~\cite[Theorem~1.1]{BouDar10}, \cite[Theorem~5.1]{ColKre98} or also~\cite[Proposition~4.4 p.~87]{TheseDarde}. We recall the identifiability result for the reader's convenience.\footnote{Note that Theorem~\ref{ThmIdentif} is true even with the following weaker assumptions: $(g,f_{b}) \in \HH^{1/2}(\partial\Omega) \times \HH^{-1/2}(\partial\Omega)$ and $\omega$ has only a continuous boundary (see~\cite[Theorem~1.1]{BouDar10}).}

\begin{theorem} \label{ThmIdentif}
The domain $\omega$ and the function $u$ that satisfy~\eqref{problemInverse} are uniquely defined by the Cauchy data $(g,f_{b}) \neq (0,0) $.
\end{theorem}

\begin{remark}
Actually we could assume that the measurement $f_{b}$ is done only on a nonempty subset~$O$ of $\partial\Omega$. All the presented result can be adapted to this case (see, e.g., \cite{Cau12}).
\end{remark}

%

In order to solve the inverse problem~\eqref{problemInverse} we will actually focus on the shape optimization problem 
\begin{equation} \label{pbOptimiz}
\displaystyle \omega^* \in \underset{\omega \in \mathcal{O}_{d_0}}{\argmin} \, J(\omega),
\end{equation}
where $J$ is the nonnegative least-square functional defined by
 \begin{equation*}
 \label{fonction:moindre-carres}
 \displaystyle  J (\omega) := \int_{\partial\Omega} \abs{ \partial_{\boldsymbol{\mathrm{n}}} u_\omega - f_b}^2 , 
 \end{equation*}
where $ u_\omega \in  \mathrm{H}^3(\priv{\Omega}{\omega})$ is the unique solution of the problem
\begin{equation} \label{pb S}
\left\lbrace 
\begin{array}{rcll}
\displaystyle - \Delta  u_{\omega} & = & 0  & \mbox{in } \, \Omega \backslash \overline{\omega}  , \\
\displaystyle  u_{\omega} &= & g  & \mbox{on } \, \partial \Omega    , \\
\displaystyle  u_{\omega} &= & 0  & \mbox{on } \, \partial \omega   .
\end{array}
\right.
\end{equation}
Indeed, the identifiability result ensures that $J(\omega)=0$ if and only if $\omega = \omega_{\rm ex}$. Finally, in order to solve numerically the shape optimization problem~\eqref{pbOptimiz}, we will now compute the shape gradient of the cost functional $J$ and apply a classical gradient descent method.

\subsection{Computation of the shape gradient.}

In order to define shape derivatives, we will use the  {\it Hadamard's method}. 
We first introduce the space of admissible deformations given by
\begin{equation} \label{AdmDeform}
\displaystyle \boldsymbol{U}:=\{\boldsymbol{V}\in\boldsymbol{\mathrm{W}}^{ 3,\infty} ;\,\mbox{\rm Supp }\boldsymbol{V}\subset\overline{\Omega_{d_{0}}}\}.
\end{equation}
In particular we are interested in the shape gradient of $J$ defined by
\begin{equation*}
\DD J(\omega) \cdot  \GV := \displaystyle \lim_{t \rightarrow 0} \frac{J\big((\GII + t \GV)(\omega)\big) - J(\omega)}{t} ,
\end{equation*}
for every $\omega \in \mathcal{O}_{d_{0}}$ and every $\GV \in \GU$. For sake of completeness, we recall the proof of the following result in Appendix~\ref{appendixproofDJ}.

\begin{proposition} \label{DJ} 
The least-square functional $J$ is differentiable at $\omega \in \mathcal{O}_{d_{0}}$ in the direction $\GV \in \GU$ with
\begin{equation} \label{formuleGradient}
   \DD J(\omega) \cdot \GV   = - \int_{\partial\omega} \partial_{\nn} u_\omega \, \partial_{\nn}w_\omega \( \GV \cdot \Gn \) ,
\end{equation}
where $w_\omega \in \HH^1(\Oomega)$ is the unique solution of the adjoint problem given by
\begin{equation} \label{PbAdjoint}
\left\{
\BA{rclll}
- \Delta w_{ \omega} & = & 0 & & \mbox{\rm in } \Oomega , \\
 w_{ \omega}  & = & 2\( \partial_{\nn} u_\omega - f_{b} \) & & \mbox{\rm on } \partial \Omega , \\
 w_{ \omega} & = & 0 & & \mbox{\rm on } \partial\omega .
\EA
\right.
\end{equation}
\end{proposition}

From the above explicit formulation of the shape gradient of $J$, we are now in a position to implement some numerical simulations based on a classical gradient descent method and we include the flip procedure introduced in this paper in order to detect in particular a multiple-component obstacle.

%
%
%


\subsection{Numerical simulations} \label{SectionNumericalSimu}

Before coming to numerical simulations, let us recall that many difficulties can be encountered in order to solve numerically Problem~\eqref{pbOptimiz}, as explained in \cite[Theorem~1]{Dambrine} (see also~\cite[Proposition~2.4]{BADR11RI}). Indeed, the gradient has not a uniform sensitivity with respect to the deformation directions. However, we use in this paper a parametric model of shape variations using piecewise B\'ezier curves which corresponds to a regularization method (as the truncated Fourier series used in~\cite{Dambrine}) allowing to overcome the ill-posedness of the inverse problem and then to solve it numerically.

{ Note that we use here piecewise B\'ezier curves that do not satisfy the $C^{2,1}$-regularity assumption made in the previous section.\footnote{ However one could retrieve the $C^{2,1}$-regularity hypothesis by imposing some additional constraints on the control points of the piecewise B\'ezier curves.} This regularity hypothesis is sufficient in order to prove rigorously the previous theoretical results. In this section dedicated to numerical simulations, we make the choice to not deal with this regularity issue since we still observe relatively good numerical reconstructions of obstacles. Besides, let us mention that the issue of knowing if the computed shape gradient (computed in particular from approximations of shapes by mesh and of solutions of PDEs by a finite element method) is an actual approximation of the genuine one is a fully-fledged question and it is not our aim to address this issue in this paper.}


\subsubsection{Framework for the numerical simulations}\label{secalgo}

The numerical simulations presented hereafter are performed in the two-dimensional case using the finite element library \textit{FreeFem++} (see~\cite{Freefem}). The exterior boundary $\partial \Omega$ is assumed to be the circle centered in the origin and of radius~$10$ and we consider the exterior Dirichlet boundary condition $ g=100$. In order to get a suitable measure $f_{b}$, we use a \textit{synthetic data}, that is, we fix a shape~$\omega_{\rm ex}$ and solve Problem~\eqref{PbExact} using a finite element method (here P2 finite element discretization) and extract the measurement $f_{b}$ by computing $\partial_{\nn} u_{\rm ex} $ on $\partial\Omega$.

Then we use a P1 finite element discretization to solve Problems~\eqref{pb S} and~\eqref{PbAdjoint} with $50$ discretization points for both the exterior boundary and each cubic B\'ezier patch describing the shape $\omega$. In order to numerically solve the optimization problem~\eqref{pbOptimiz}, we use the following classical gradient descent algorithm and we include the flip procedure at Step~\eqref{step3}. \\
~\\
\textbf{Algorithm $\mathcal{A}$}
\begin{enumerate}
\item Fix $k=0$, fix an initial shape $\omega_{0}$, fix a maximal number $M \in \N^*$ of iterations and fix $\lambda \geq 1$ a given tolerance coefficient for the flip procedure (see Step~\eqref{flipoupasflip}, $\lambda$ should be chosen close to~$1$).
\item Control the size of the patches of $\omega_{k}$ (see Remark~\ref{remcontrolsize}). \label{step2}
\item Scan $\omega_{k}$ looking for intersecting control polygons (see Section~\ref{intersectdetection}): \label{step3} \vspace{4pt}
\begin{enumerate}
\item in the case of no intersecting control polygons, go to Step \eqref{step4}; \vspace{4pt}
\item in the case of intersecting control polygons: \vspace{4pt}
\begin{enumerate}
\item apply the flip procedure and obtain a multiple-component shape $\omega^1_k \cup \omega^2_k$; \\[2pt]
recall that $\omega^1_k \cup \omega^2_k$ is not self-intersecting from Assumption (A3);\vspace{4pt}
\item compute $J(\omega^1_k \cup \omega^2_k)$ and $J(\omega_k)$, { and set $J(\omega_k)=+\infty$ if $\omega_k$ is self-intersecting}: \label{flipoupasflip} \vspace{4pt}
\begin{enumerate}
\item if $J(\omega^1_k \cup \omega^2_k) < \lambda J(\omega_k) $, do $\omega_k \leftarrow \omega^1_k \cup \omega^2_k$; \vspace{4pt}
\item else, go to Step \eqref{step4}. \label{flipoupasflip2}
\end{enumerate}
\end{enumerate}
\end{enumerate}
\item Solve Problems \eqref{pb S} and~\eqref{PbAdjoint} with $\omega = \omega_{k}$. \label{step4}
\item Compute the shape gradient $\DD J (\omega_{k})$ from Formula~\eqref{formuleGradient}.
\item Move the control points of the shape, that is, do $\omega_{k+1} \leftarrow \omega_{k} - \alpha_k \DD J (\omega_{k})$, where $\alpha_k$ is a small positive coefficient chosen, e.g., by a classical line search.
\item Do $k \leftarrow k+1$ and get back to Step~\eqref{step2} while $k < M$.
\end{enumerate}





\subsubsection{First simulations: detection of smooth and convex shapes}

\begin{figure}[h]
  \centering
  \subfigure[Detection of a circle]{
  \includegraphics[trim={0 0 0 0cm},scale=0.4]{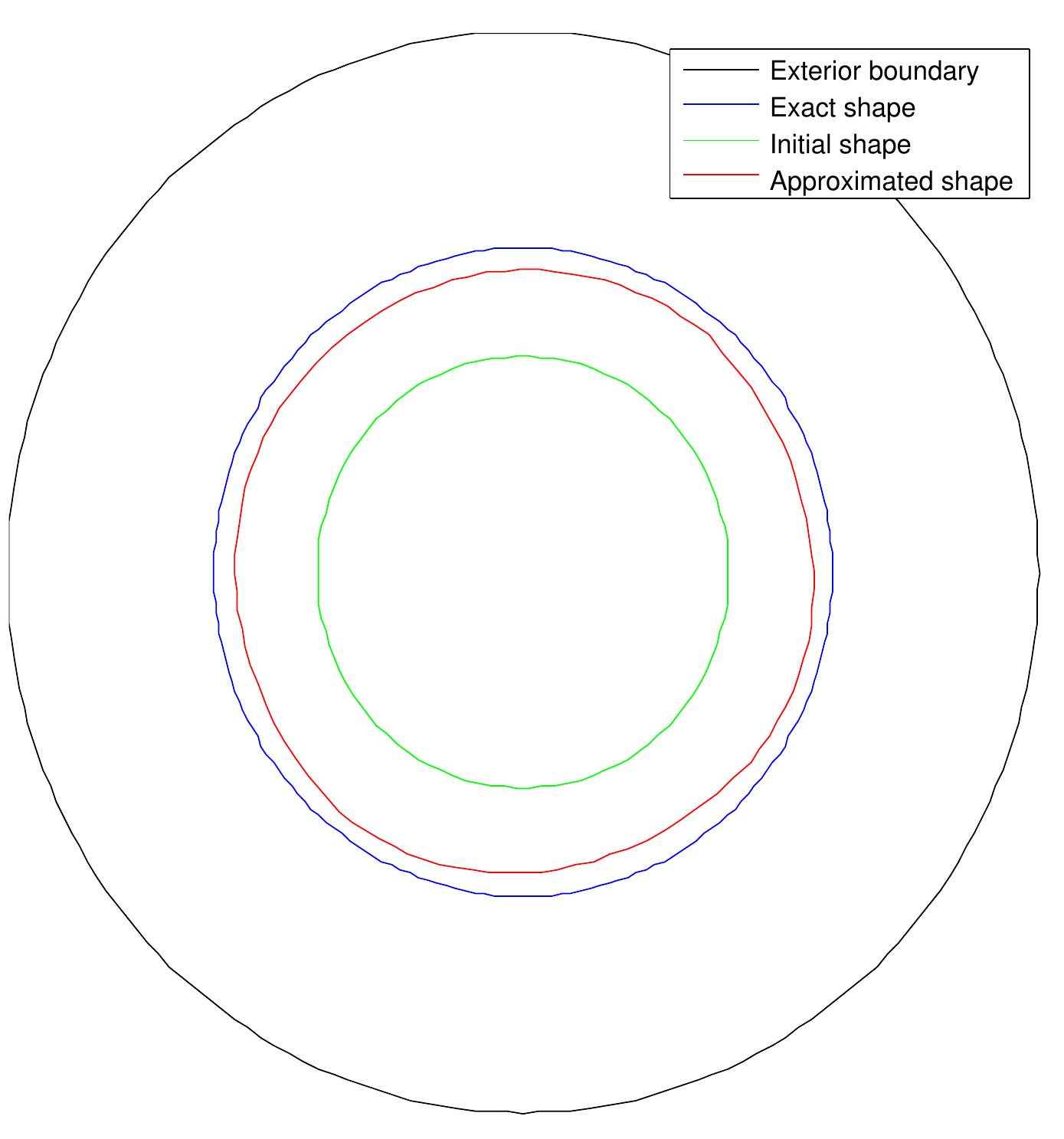}
  \label{simuCercle}}
\quad
  \subfigure[Detection of an ellipse]{
  \includegraphics[trim={0 0 0 0cm},scale=0.4]{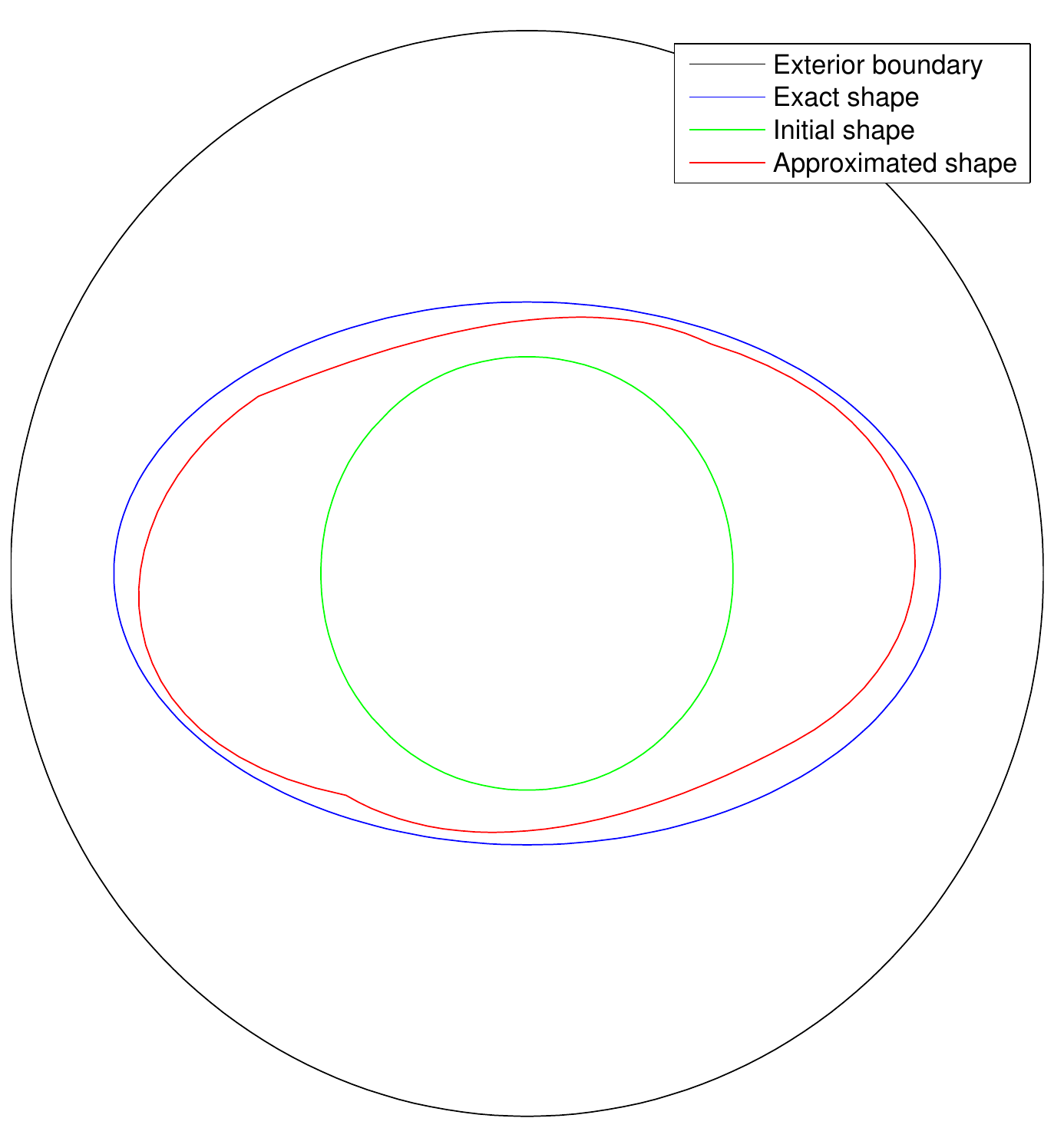}
  \label{simuEllipse}}
\caption{Detection of convex and smooth obstacles.}
\label{simuCercleEllipse}
\end{figure}

We first test Algorithm $\mathcal{A}$ on the problem of detecting one smooth convex object. Precisely, we begin by detecting the circle centered at the origin and of radius $6$ and the ellipse $\{ (8\cos \theta,\ 5\sin \theta), \theta \in [0,2\pi]\}$ using four cubic B\'ezier patches. Numerical simulations are performed and depicted in Figure~\ref{simuCercleEllipse}.

\subsubsection{Detection of a non-smooth shape and of a non-convex shape}

\begin{figure}[h]
  \centering
  \subfigure[Detection of a square]{
  \includegraphics[trim={0 0 0 0cm},scale=0.4]{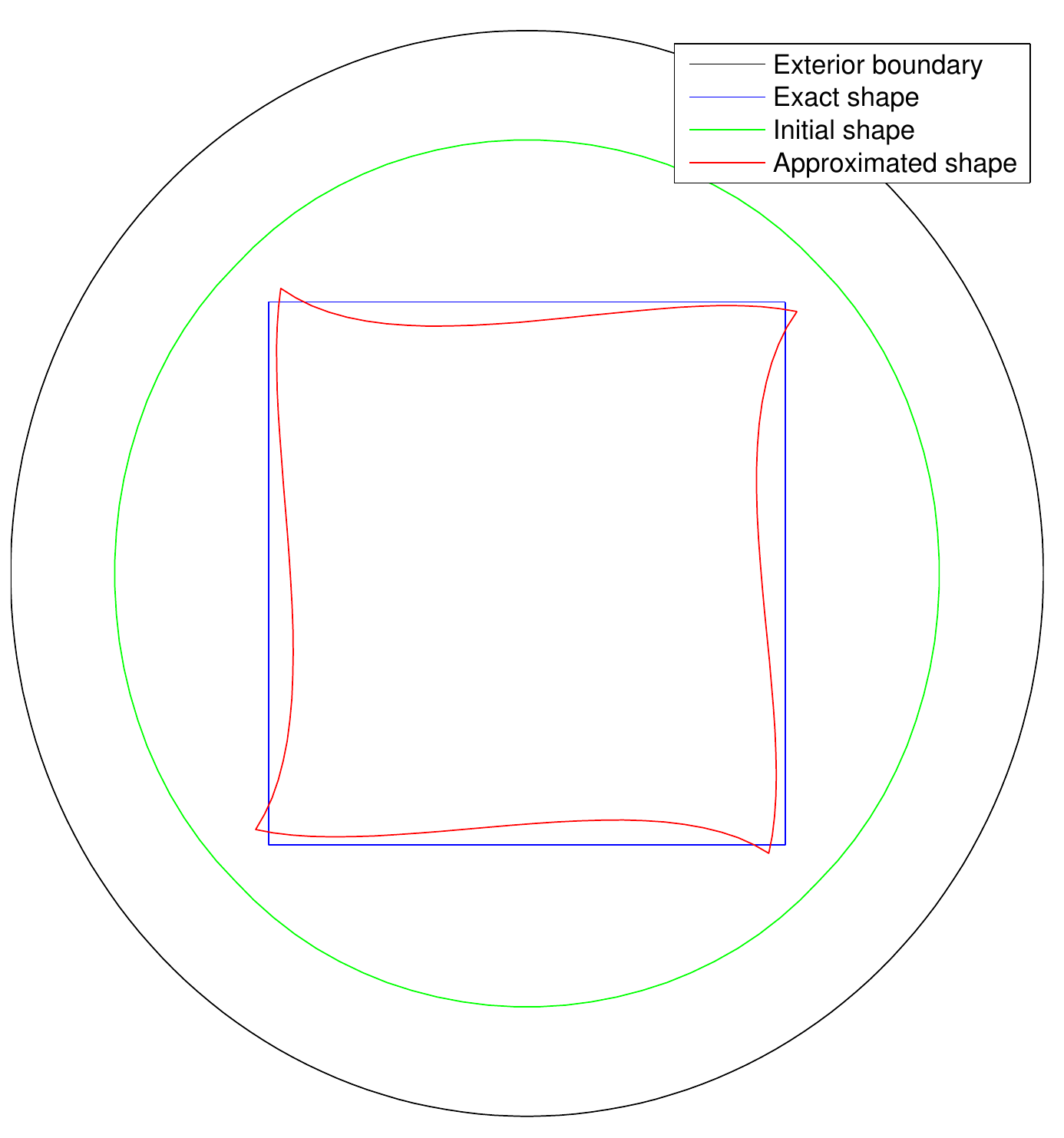}
  \label{simuCarre}}
\quad
  \subfigure[Detection of a non-convex shape.]{
  \includegraphics[trim={0 0 0 0cm},scale=0.4]{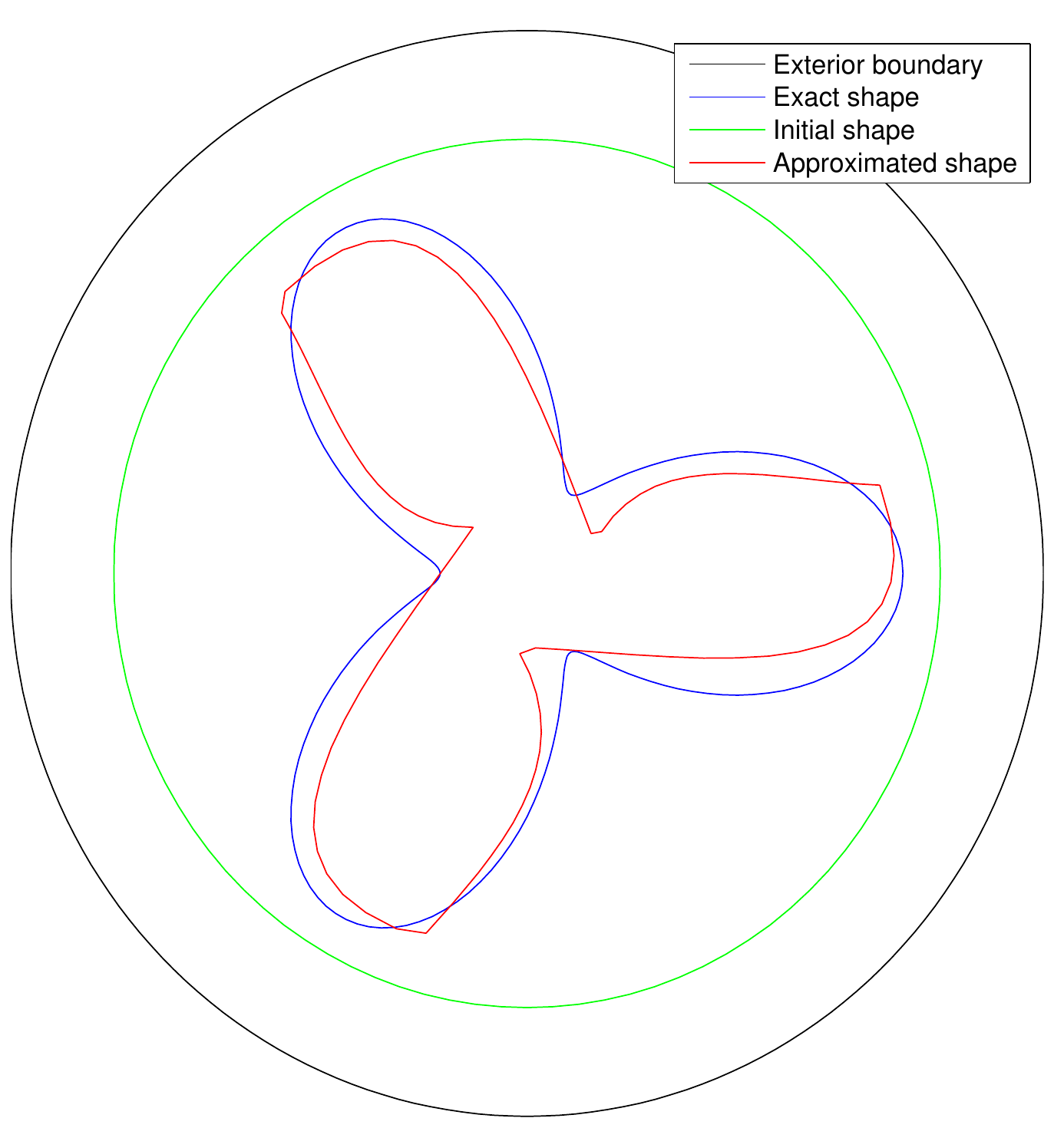}
  \label{simuFormey}}
\caption{Detection of a non-smooth obstacle and of a non-convex obstacle.}
\label{simunonsmoothnonconvex}
\end{figure}

\begin{figure}[h]
  \centering
  \includegraphics[scale=0.4]{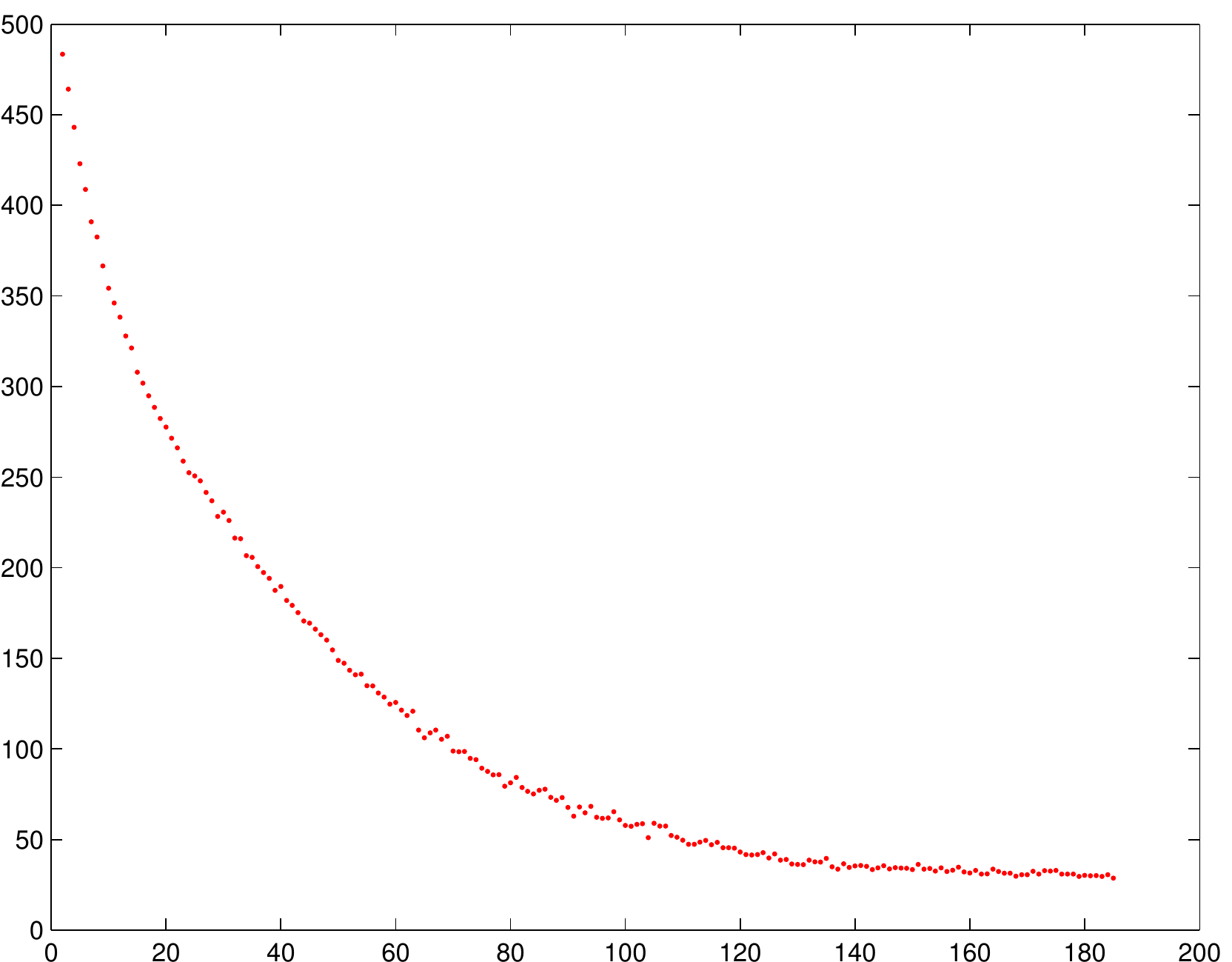}
  \caption{Evolution of the objective function for the detection of the square.}
  \label{evolution_critere_carre}
\end{figure}

We test now Algorithm $\mathcal{A}$ on the problem of detecting a non-smooth shape and of detecting a non-convex shape (see Figure~\ref{simunonsmoothnonconvex}). Precisely we first consider the square of side $10$ and centered at the origin and we use four cubic B\'ezier patches. As one can see in Figure~\ref{simuCarre}, each B\'ezier patch detects a side of the square. Figure~\ref{evolution_critere_carre} shows the decrease of the objective function during the simulation. Secondly, in Figure~\ref{simuFormey}, we consider the non-convex shape parameterized by $\{ (2.8(1.6+ \cos(3\theta)) \cos(\theta), 2.8(1.6+ \cos(3\theta)) \cos(\theta) ), \theta\in[0,2\pi] \}$, using six cubic B\'ezier patches.\footnote{This shape is also considered in~\cite[Figure~4]{CaubetDambrineKateb} where authors obtained the convex hull of the shape. However, note that the authors used a different method where the descent direction is obtained by solving a boundary value problem involving the kernel of the shape gradient.} 



\subsubsection{Detection of a two-component obstacle starting from a one-component shape}\label{sectionseveralobstacles}
\begin{figure}[h]
  \centering
\subfigure[Initial shape]
  {
    \includegraphics[width=0.31\textwidth,scale=0.35]{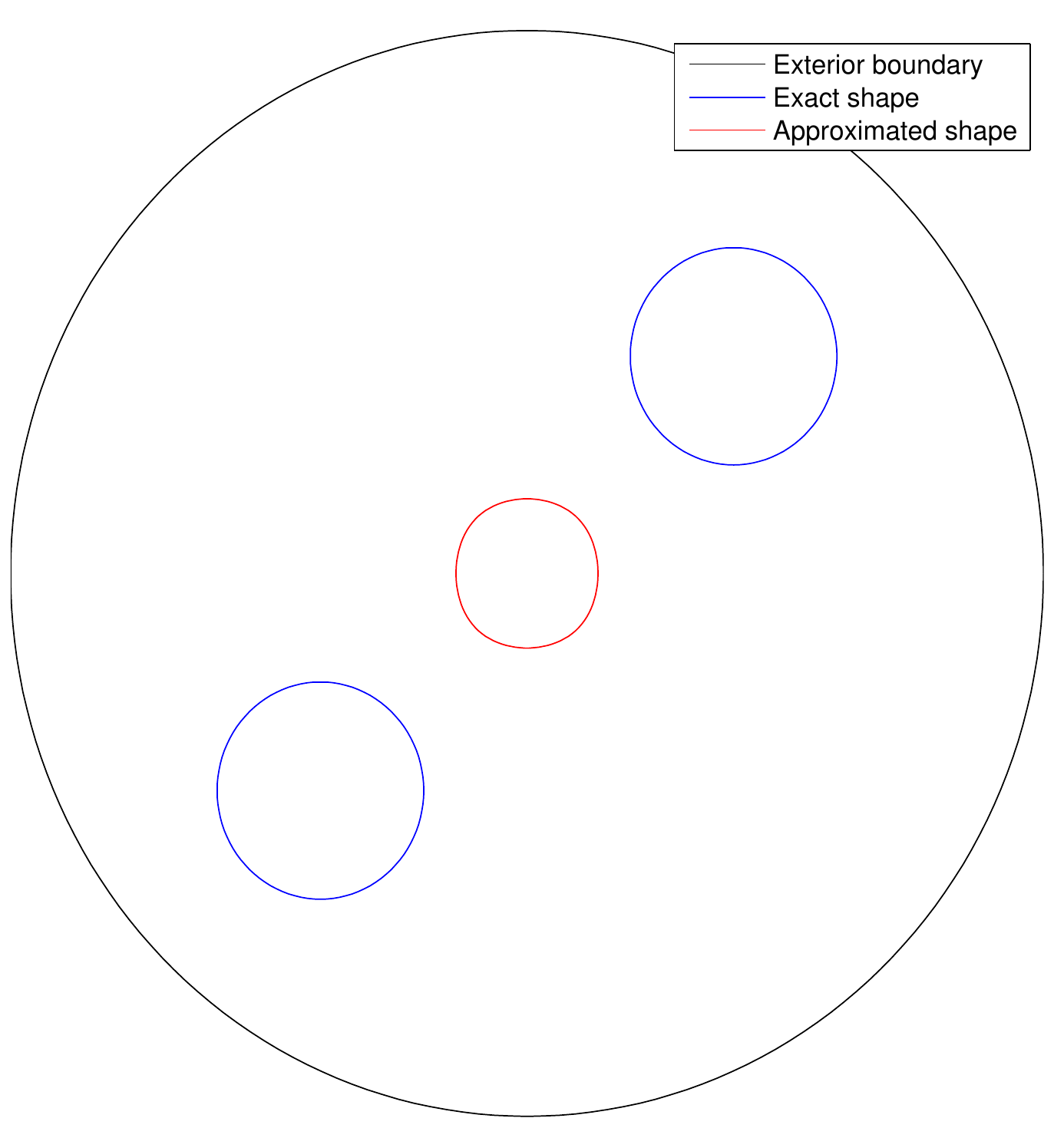}
    \label{simuFlip1}
  }
\quad
\subfigure[Intersecting control polygons]
  {
    \includegraphics[width=0.31\textwidth,scale=0.35]{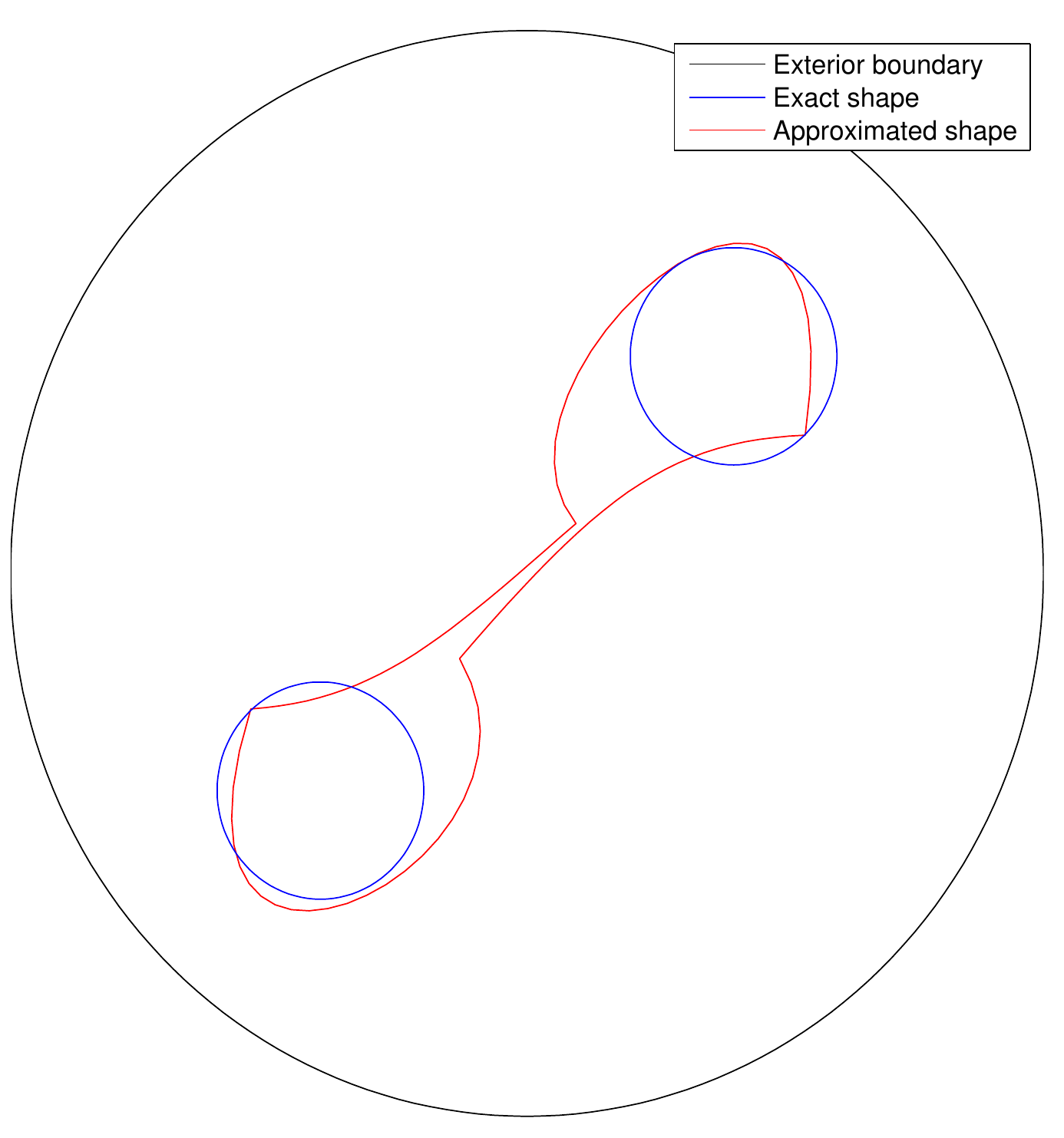}
    \label{simuFlip2}
  }
\subfigure[Flip procedure]
  {
    \includegraphics[width=0.31\textwidth,scale=0.35]{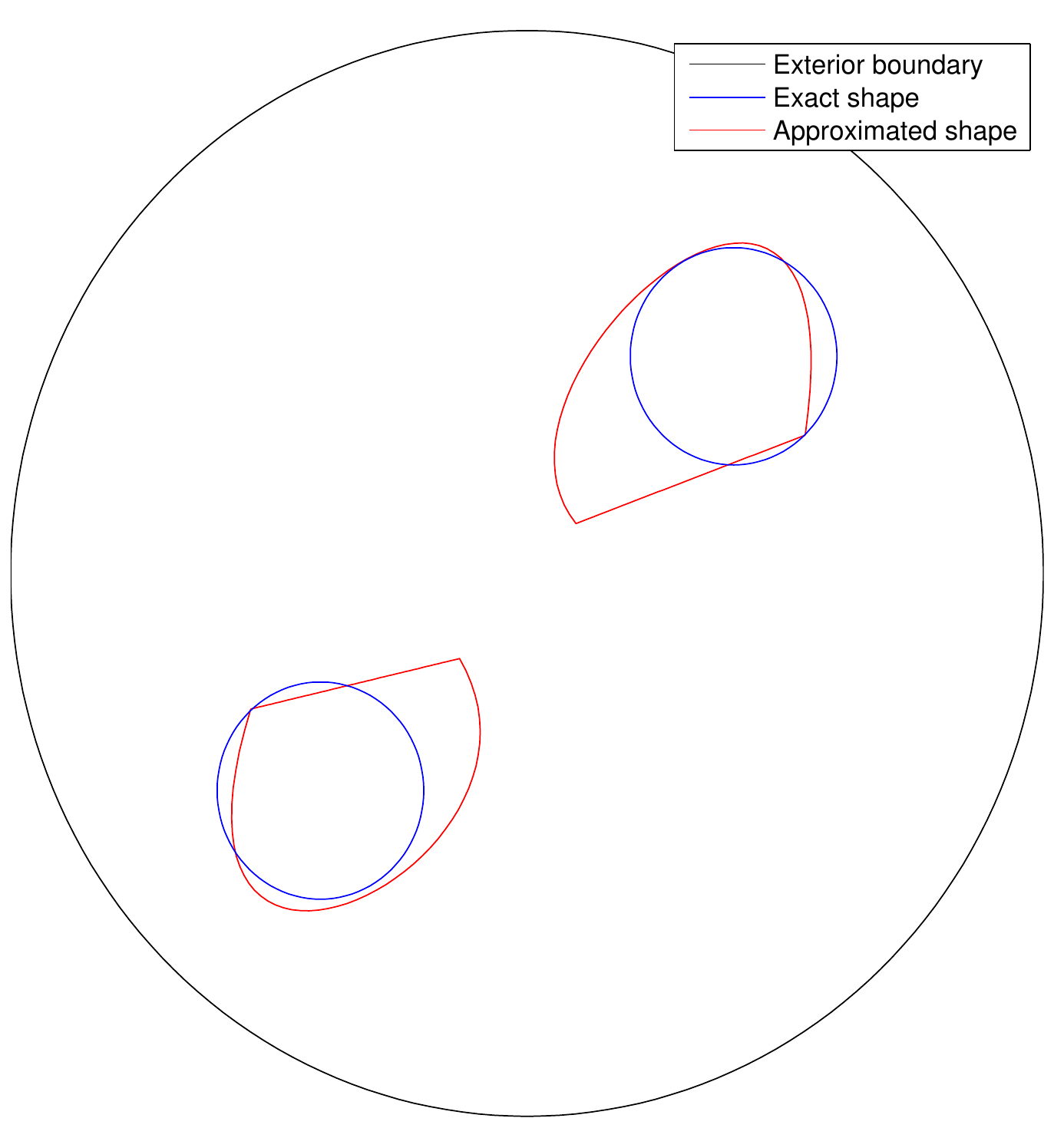}
    \label{simuFlip3}
  }
\quad
\subfigure[Final shape]
  {
    \includegraphics[width=0.31\textwidth,scale=0.35]{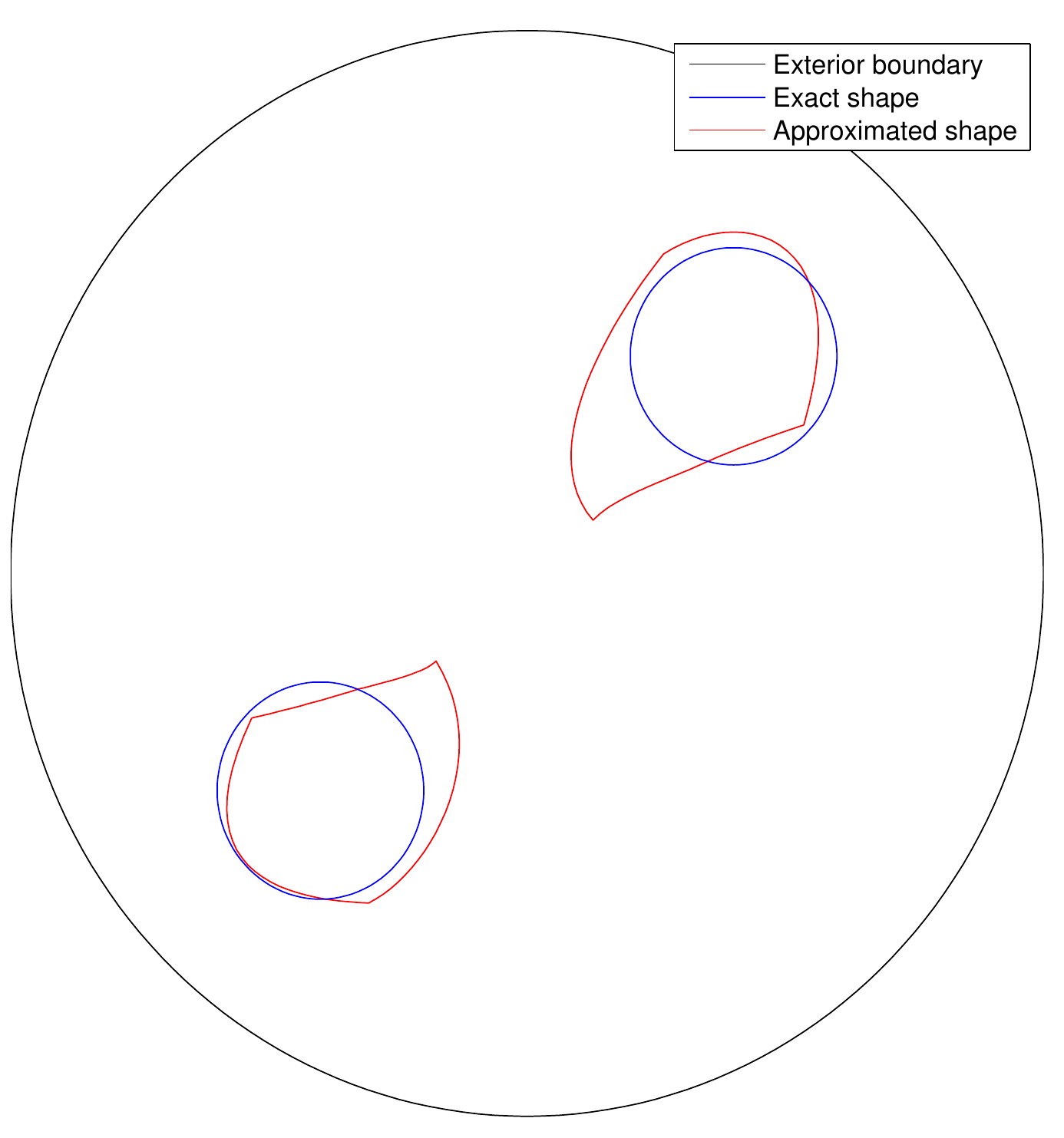}
    \label{simuFlip4}
  }
\caption{Detection of two obstacles starting from a one-component shape}
\label{simuFlip}
\end{figure}

In this section we test the flip procedure introduced in Section~\ref{SectionFlip} in order to detect a two-component shape starting from a one-component initial shape. We consider two circles of radius $2$ centered at $(-4,-4)$ and $(4,4)$. We present different states of the algorithm in Figure~\ref{simuFlip}. The initial B\'ezier shape consists in a single component with four cubic B\'ezier patches, located at the center (Figure~\ref{simuFlip1}). The shape grows and surrounds the two objects until two control polygons intersect each other (Figure~\ref{simuFlip2}). The flip procedure is performed and the shape is divided in two connected components (Figure~\ref{simuFlip3}). At the end, the algorithm provides an approximation of the two obstacles (Figure~\ref{simuFlip4}). 

\begin{figure}[h]
  \centering
  \includegraphics[scale=0.5]{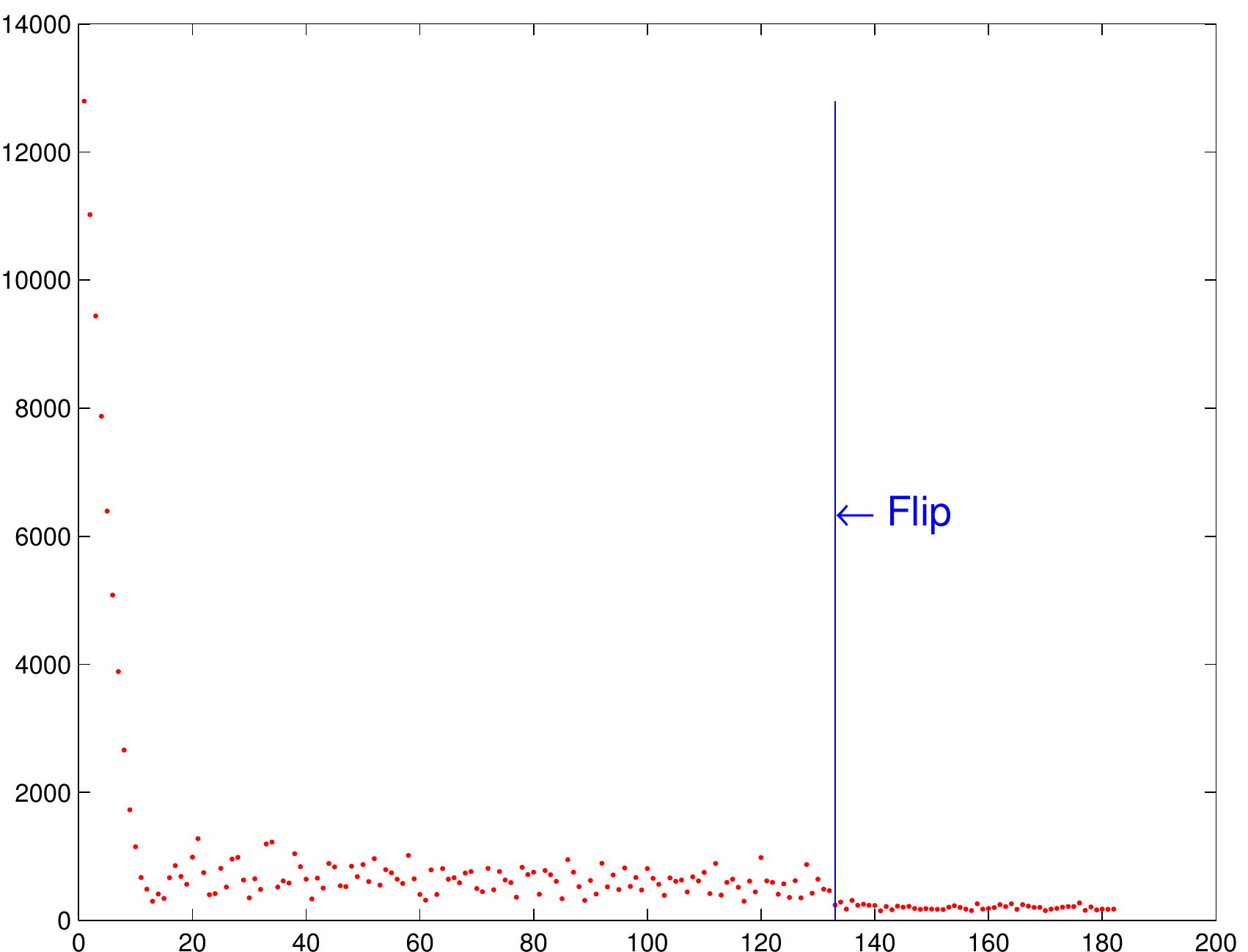}
  \caption{Evolution of the objective function for the detection of a two-component shape.}
  \label{evolution_critere_flip}
\end{figure}

Figure~\ref{evolution_critere_flip} depicts the evolution of the objective function during this simulation. 
One can note a change of behavior after Iteration $133$ which corresponds to the performance of the flip procedure. Precisely, the algorithm finds in a first place a local minimizer at Iteration $13$, which corresponds to a one-component minimizer. After oscillations around this local minimum, the flip procedure is performed and the functional decreases and stabilizes around a two-component minimizer.



\subsubsection{Checking the objective function value after a flip procedure}\label{sectionflipoupasflip}
\begin{figure}[h]
  \centering
  \subfigure[Before the flip procedure, $J(\omega_k) = 3211$]{
  \includegraphics[trim={0 0 0 0cm},scale=0.4]{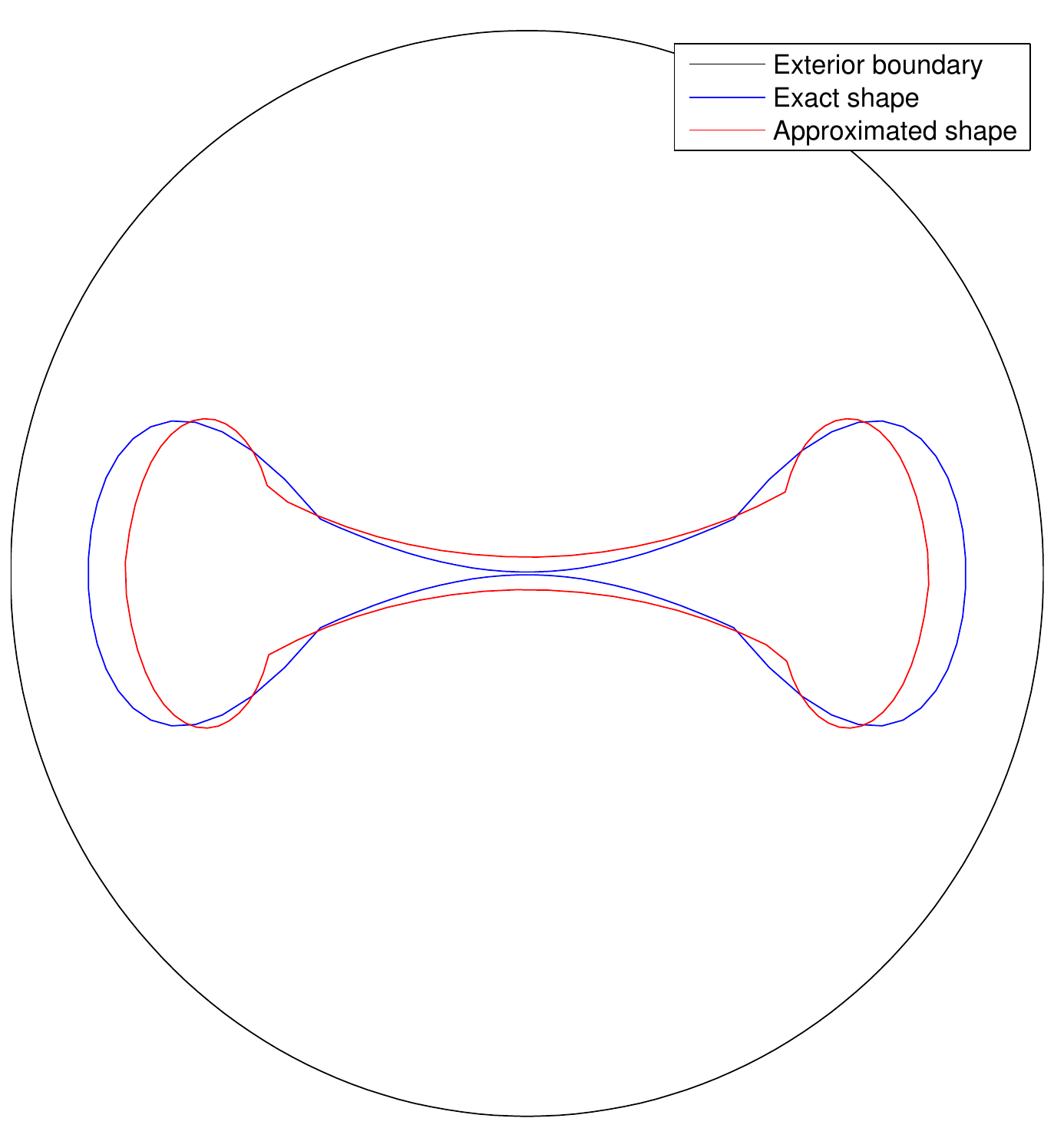}
  \label{simuSautAvantFlip}}
\quad
  \subfigure[After the flip procedure, $J(\omega^1_k \cup \omega^2_k) = 3579$]{
  \includegraphics[trim={0 0 0 0cm},scale=0.4]{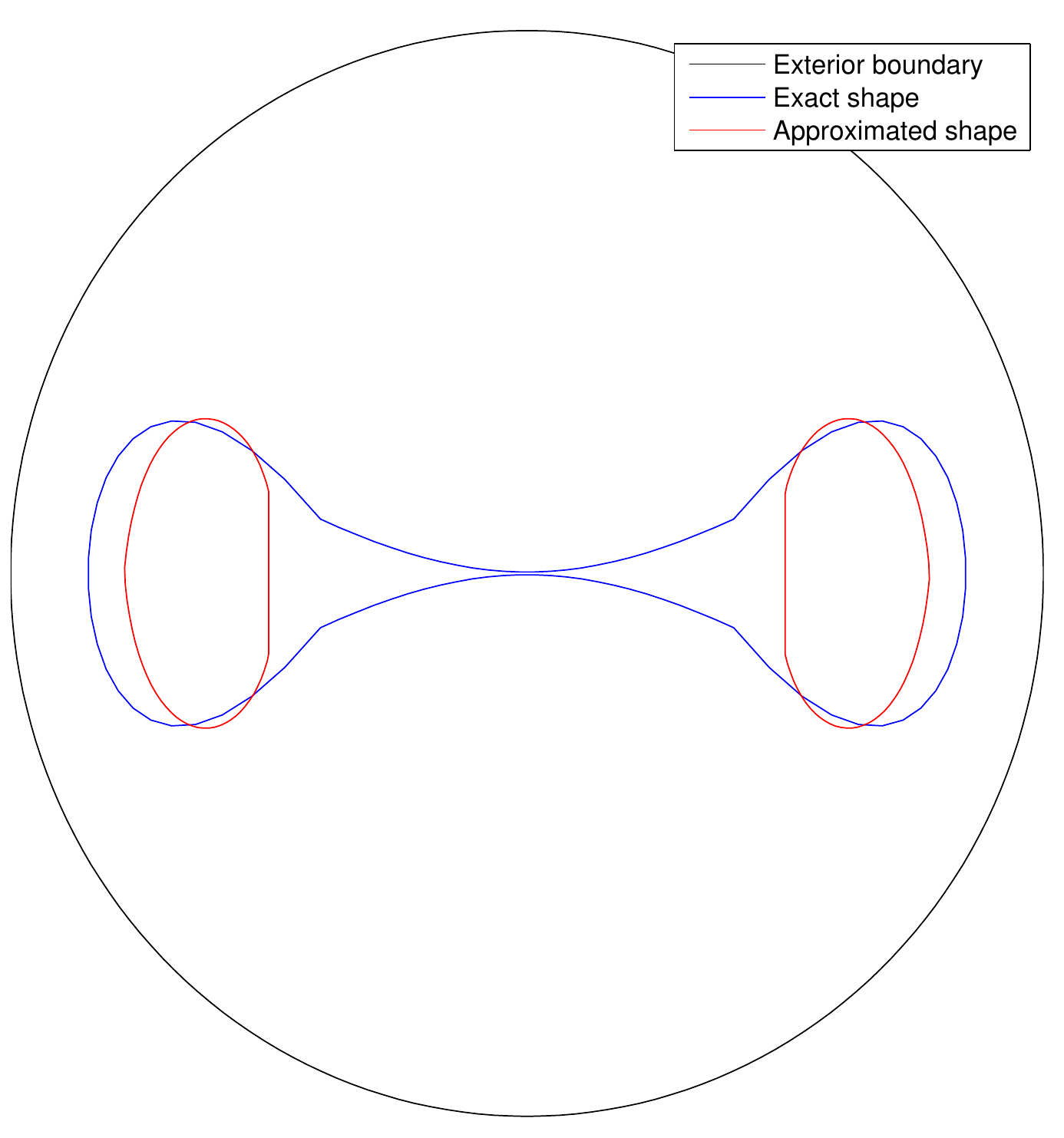}
  \label{simuSautApresFlip}}
\caption{The objective function value significantly increases whenever the flip procedure is not needed.}
\label{cancelFlip}
\end{figure}

In Algorithm $\mathcal{A}$, Step~\eqref{flipoupasflip} makes sure that, whenever a flip procedure is performed, the objective function value does not significantly increase. If $J(\omega^1_k \cup \omega^2_k) \geq \lambda J(\omega_k) $ (for instance $\lambda = 1.1$), then we consider that adding another component to the shape is not a wise choice and we cancel the flip procedure. This situation typically occurs when the target shape has a single component with a very thin part (\textit{i.e.} two parts of its boundary are very close to each other). In such a case Algorithm $\mathcal{A}$ probably leads to two control polygons intersecting each other and to a flip performance, while the target shape has a single component. We present an example of such a situation in Figure~\ref{cancelFlip}. The obstacle is composed of one component with a very thin part and the current shape $\omega_k$ of the algorithm has two control polygons intersecting each other. The objective function value before the flip procedure is $J(\omega_k) = 3211$ and after the flip procedure, it has increased to $J(\omega^1_k \cup \omega^2_k) = 3579$. Since the ratio is greater than~$\lambda$, the algorithm cancels the flip procedure and goes to Step~\eqref{step4}.

\begin{remark}
{ In Algorithm $\mathcal{A}$, if the previous situation $\omega_k$ is restored after the performance of a flip procedure (that is, in the case~\ref{flipoupasflip2}), then the obstacle is highly likely composed of one component with a very thin part. Next, the gradient flow makes evolve the approximated shape $\omega_k$ with a small deformation step to a better approximation $\omega_{k+1}$ of the target shape. This would lead to a new perfomance of a flip procedure on $\omega_{k+1}$. Actually, in such a situation, a flip procedure is performed/cancelled at each iteration until the end. }
\end{remark}

\section{Conclusion and perspectives}

In this paper we studied the use of a piecewise B\'ezier parameterization for the representation of two-dimensional shapes in geometric approximation based on successive shape deformations. We proposed procedures in order to manipulate this parameterization and showed how to manage changes of topology and so multiple-component shape approximations. We applied this approach to a problem of multiple-inclusion detection and performed numerical simulations using \textit{FreeFem++}. The computational efficiency of the method arises from the simplicity and the flexibility of the proposed parameterization. 

We considered in this paper a two-dimensional problem, but the challenging extension to the three-dimensional case may be interesting and the algorithmic contents could be generalized. As a conclusion, the implementation described here was made for an experimental purpose and a complete and optimized implementation may be proposed.

{ 
\paragraph{Acknowledgment.}
 
The authors would like to thank the anonymous reviewers for their valuable comments and suggestions in order to improve the quality of the paper. }

\appendix

\section{Proof of Proposition~\ref{DJ}}\label{appendixproofDJ}

We detail here the classical proof of Proposition~\ref{DJ} for the reader's convenience. For any $\GV \in \GU$ (where~$\GU$ is defined by~\eqref{AdmDeform}), we introduce the perturbed domain $\omega_{t} := (\GII + t \GV ) (\omega)$ and the functional~$j$ defined for all $t\in[0,T)$ by $j(t):=J(\omega_{t})$ and we consider the unique solution $u_{t} \in  \HH^3(\priv{\Omega}{\omega_{t}})$ of the perturbed problem
$$
\left\{
\BA{rclll}
- \Delta u_t & = & 0 & & \mbox{\rm in } \Oomegat , \\
 u_t & = &  g & & \mbox{\rm on } \partial \Omega,  \\
 u_t & = & 0 & & \mbox{\rm on } \partial\omega_{t} .

\EA
\right.
$$
Let us recall the definition of the shape derivative in our situation (see~\cite{HP} for details). We introduce
$$
\boldsymbol{\mathcal{U}}:= \left\lbrace \Gtheta \in \GU  ; \,   \Vert \Gtheta \Vert_{ 3,\infty}    < \min \( \frac{d_{0}}{3} , 1 \)  \right\rbrace   
$$
and, for any $\Gtheta \in \boldsymbol{\mathcal{U}}$, we consider the unique solution $u_{\theta} \in \HH^3(\priv{\Omega}{\omega_{\theta}})$ of the perturbed problem
$$
\left\{
\BA{rclll}
- \Delta u_{\theta} & = & 0 & & \mbox{\rm in } \priv{\Omega}{\omega_{\theta}} , \\
 u_{\theta} & = & g & & \mbox{\rm on } \partial \Omega,  \\
 u_{\theta} & = & 0 & & \mbox{\rm on } \partial\omega_{\theta} ,
\EA
\right.
$$
where $\omega_{\theta} := (\GII + \Gtheta ) (\omega)$. Then,
\begin{itemize}
\item if the mapping $\Gtheta \in \boldsymbol{\mathcal{U}} \mapsto u_{\theta} \circ (\GII + \Gtheta)  \in \HH^1(\Oomega) $ is Fr\'echet differentiable at $\0$, we say that $\Gtheta \mapsto u_{\theta} $ possesses a \textit{total first variation} (or derivative) at $\0$. In such a case, this total first derivative at $\0$ in the direction $\Gtheta$ is denoted by $ \stackrel{.}{u}_{\theta} $ and is called \textit{material derivative} (or \textit{Lagrangian derivative)};
\item if, for every $\mathscr{D} \subset \subset \Oomega$, the mapping $\Gtheta \in \boldsymbol{\mathcal{U}} \mapsto \restriction{u_{\theta}}{\mathscr{D}}  \in \HH^1(\mathscr{D}) $ is Fr\'echet differentiable at $\0$, we say that $\Gtheta \mapsto  u_{\theta} $ possesses a \textit{local first variation} (or derivative) at $\0$. In such a case, this local first derivative at $\0$ in the direction~$\Gtheta$ is denoted by $u'_{\theta} $, is called \textit{shape derivative} (or \textit{Eulerian derivative}) and is well defined in the whole domain $\Oomega$:
\BEAN
u'_{\theta}  = \frac{d}{d t} \restriction{\( \restriction{u_{t \Gtheta}}{\mathscr{D}} \)}{t=0} & & \mbox{in each } \, \mathscr{D} \subset \subset \Oomega.
\EEAN
\end{itemize}
In the sequel, let $\GV \in \boldsymbol{\mathcal{U}}$ and let $u'$ be the local first variation $u'_{V} $ which is referred as the shape derivative of the state.

The differentiability of the cost functional $J$ is directly obtained from the existence of the shape derivative of the state $u$ given for example in~\cite[Theorem~5.3.1]{HP}. Notice that in~\cite[Theorem~5.3.1]{HP}, the result claims the differentiability of $t\in [0,T] \mapsto \tilde u_{t} \in \LL^2(\Omega)$, where $\tilde u_{t}$ is an extension of $u_{t}$ in~$\Omega$. Since we want to obtain the differentiability of $t\in [0,T] \mapsto \tilde u_{t} \in  \HH^2(\Omega)$ (in order to differentiate properly the functional~$J$), we have here to work with the mentioned spaces, that is with domains with a~$C^{2,1}$ boundary (and not only Lipschitz) and perturbations~$\GV$ which belong to~$ \GWW^{3,\infty} (\R^2)$ (and not only to~$\GWW^{1,\infty} (\R^2)$).

Moreover we can easily characterize the shape derivative $u'\in\HH^1(\Oomega)$ as the solution of the following problem (see again for example~\cite[Theorem~5.3.1]{HP}):
\BE \label{PbDeriv}
\left\{
\BA{rclll}
- \Delta u' & = & 0 & & \mbox{\rm in } \Oomega , \\
 u' & = & 0 & & \mbox{\rm on } \partial \Omega , \\
 u' & = & - \partial_{\nn} u \( \GV \cdot \Gn \) & & \mbox{\rm on } \partial\omega .
\EA
\right.
\EE

Then by differentiation under the sum sign, we obtain
$$
j'(0) = 2 \int_{\partial\Omega} \partial_{\nn} u' (\partial_{\nn} u - f_{b}) .
$$
Using the weak formulation of Problem~\eqref{PbDeriv} solved by $u'$ with $w$ as a test function, we obtain
\begin{equation*}
\int_{\Oomega} \nabla u' \cdot \nabla w - \int_{\partial(\Oomega)} w \, \partial_{\nn} u'  = 0
\end{equation*}
and using the weak formulation of the adjoint Problem~\eqref{PbAdjoint} solved by $w$ with $u'$ as a test function, we obtain
\begin{equation*}
\int_{\Oomega} \nabla w \cdot \nabla u' - \int_{\partial(\Oomega)} u' \, \partial_{\nn} w  = 0 .
\end{equation*}
Finally, using the boundary conditions, the proof is complete.

\section{Detection of one obstacle starting from a two-component shape}\label{2for1}
In this paper we have introduced the flip procedure as a method that enables to divide a one-component shape into a two-component shape. Actually the flip procedure can be easily adapted in order to perform the reverse operation, that is, to merge a two-component shape into a one-component shape (see Figure~\ref{merge}). 

We focus now on the detection of the one-component shape $\{ (4\cos \theta , 6+2.5\sin \theta), \theta\in[0,2\pi] \}$ and we start Algorithm $\mathcal{A}$ with a two-component shape. 
We present different states of the algorithm in Figure~\ref{simuFlipMerge}. At the end, the algorithm provides an approximation of the one-component obstacle.

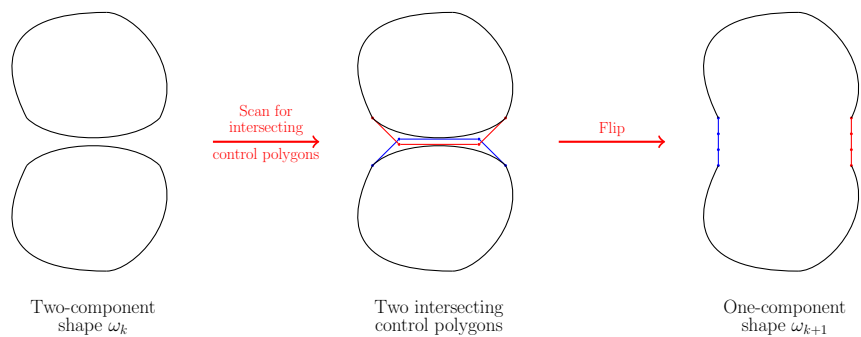
\begin{figure}[h]
\centering
\begin{tikzpicture}[scale=0.35, every node/.style={scale=0.35}]
\node (P1) at (2.000000,4.000000) {};
\node (P2) at (1.000000,2.000000) {};
\node (P3) at (1.000000,0.000000) {};
\node (P4) at (5.000000,0.000000) {};
\node (P5) at (6.000000,0.000000) {};
\node (P6) at (8.000000,2.000000) {};
\node (P7) at (7.000000,4.000000) {};
\node (P8) at (6.000000,5.000000) {};
\node (P9) at (3.000000,5.000000) {};
\node (Q1) at (2.000000,5.800000) {};
\node (Q2) at (1.000000,7.800000) {};
\node (Q3) at (1.000000,9.800000) {};
\node (Q4) at (5.000000,9.800000) {};
\node (Q5) at (6.000000,9.800000) {};
\node (Q6) at (8.000000,7.800000) {};
\node (Q7) at (7.000000,5.800000) {};
\node (Q8) at (6.000000,4.800000) {};
\node (Q9) at (3.000000,4.800000) {};
\draw[black] (2.000000,4.000000) .. controls (1.000000,2.000000) and (1.000000,0.000000) .. (5.000000,0.000000);
\draw[black] (5.000000,0.000000) .. controls (6.000000,0.000000) and (8.000000,2.000000) .. (7.000000,4.000000);
\draw[black] (7.000000,4.000000) .. controls (6.000000,5.000000) and (3.000000,5.000000) .. (2.000000,4.000000);
\draw[black] (2.000000,5.800000) .. controls (1.000000,7.800000) and (1.000000,9.800000) .. (5.000000,9.800000);
\draw[black] (5.000000,9.800000) .. controls (6.000000,9.800000) and (8.000000,7.800000) .. (7.000000,5.800000);
\draw[black] (7.000000,5.800000) .. controls (6.000000,4.800000) and (3.000000,4.800000) .. (2.000000,5.800000);
\node (x) at (4.5,-1.4) {\huge Two-component};
\node (x) at (4.5,-2.1) {\huge shape $\omega_k$};

\node(flip) at (11,6.1) {\textcolor{red}{\LARGE Scan for}};
\node (flip2) at (11,5.4) {\textcolor{red}{\LARGE intersecting}};
\node (flip3) at (11,4.4) {\textcolor{red}{\LARGE control polygons}};
\draw[red,thick,->] (9,4.9) -- (13,4.9);

\node[circle,fill=blue,scale=0.3] (P1) at (13+2.000000,4.000000) {};
\node (P2) at (13+1.000000,2.000000) {};
\node (P3) at (13+1.000000,0.000000) {};
\node (P4) at (13+5.000000,0.000000) {};
\node (P5) at (13+6.000000,0.000000) {};
\node (P6) at (13+8.000000,2.000000) {};
\node[circle,fill=blue,scale=0.3] (P7) at (13+7.000000,4.000000) {};
\node[circle,fill=blue,scale=0.3] (P8) at (13+6.000000,5.000000) {};
\node[circle,fill=blue,scale=0.3] (P9) at (13+3.000000,5.000000) {};
\node[circle,fill=red,scale=0.3] (Q1) at (13+2.000000,5.800000) {};
\node (Q2) at (13+1.000000,7.800000) {};
\node (Q3) at (13+1.000000,9.800000) {};
\node (Q4) at (13+5.000000,9.800000) {};
\node (Q5) at (13+6.000000,9.800000) {};
\node (Q6) at (13+8.000000,7.800000) {};
\node[circle,fill=red,scale=0.3] (Q7) at (13+7.000000,5.800000) {};
\node[circle,fill=red,scale=0.3] (Q8) at (13+6.000000,4.800000) {};
\node[circle,fill=red,scale=0.3] (Q9) at (13+3.000000,4.800000) {};
\draw[blue] (P7) -- (P8) -- (P9) -- (P1);
\draw[red] (Q7) -- (Q8) -- (Q9) -- (Q1);
\draw[black] (13+2.000000,4.000000) .. controls (13+1.000000,2.000000) and (13+1.000000,0.000000) .. (13+5.000000,0.000000);
\draw[black] (13+5.000000,0.000000) .. controls (13+6.000000,0.000000) and (13+8.000000,2.000000) .. (13+7.000000,4.000000);
\draw[black] (13+7.000000,4.000000) .. controls (13+6.000000,5.000000) and (13+3.000000,5.000000) .. (13+2.000000,4.000000);
\draw[black] (13+2.000000,5.800000) .. controls (13+1.000000,7.800000) and (13+1.000000,9.800000) .. (13+5.000000,9.800000);
\draw[black] (13+5.000000,9.800000) .. controls (13+6.000000,9.800000) and (13+8.000000,7.800000) .. (13+7.000000,5.800000);
\draw[black] (13+7.000000,5.800000) .. controls (13+6.000000,4.800000) and (13+3.000000,4.800000) .. (13+2.000000,5.800000);
\node (x) at (13+4.5,-1.4) {\huge Two intersecting};
\node (x) at (13+4.5,-2.1) {\huge control polygons};

\node (flip1) at (24,5.4) {\textcolor{red}{\LARGE Flip}};
\draw[red,thick,->] (13+9,4.9) -- (13+13,4.9);

\node[circle,fill=blue,scale=0.3] (P1) at (13+13+2.000000,4.000000) {};
\node (P2) at (13+13+1.000000,2.000000) {};
\node (P3) at (13+13+1.000000,0.000000) {};
\node (P4) at (13+13+5.000000,0.000000) {};
\node (P5) at (13+13+6.000000,0.000000) {};
\node (P6) at (13+13+8.000000,2.000000) {};
\node (P7) at (13+13+8.000000,7.800000) {};
\node[circle,fill=red,scale=0.3] (P8) at (13+13+7.000000,4.000000) {};
\node[circle,fill=red,scale=0.3] (P9) at (13+13+7.000000, 4 + 0.6) {};
\node[circle,fill=red,scale=0.3] (P10) at (13+13+7.000000, 4 + 1.2) {};
\node[circle,fill=red,scale=0.3] (P11) at (13+13+7.000000,5.800000) {};
\node (P12) at (13+13+6.000000,9.800000) {};
\node (P13) at (13+13+5.000000,9.800000) {};
\node (P14) at (13+13+1.000000,9.800000) {};
\node (P15) at (13+13+1.000000,7.800000) {};
\node[circle,fill=blue,scale=0.3] (P16) at (13+13+2.000000,5.8 - 0.6) {};
\node[circle,fill=blue,scale=0.3] (P17) at (13+13+2.000000,5.8 - 1.2) {};
\node[circle,fill=blue,scale=0.3] (P18) at (13+13+2.000000,5.800000) {};
\draw[black] (13+13+2.000000,4.000000) .. controls (13+13+1.000000,2.000000) and (13+13+1.000000,0.000000) .. (13+13+5.000000,0.000000);
\draw[black] (13+13+5.000000,0.000000) .. controls (13+13+6.000000,0.000000) and (13+13+8.000000,2.000000) .. (13+13+7.000000,4.000000);
\draw[red] (13+13+7,4) -- (13+13+7,5.8);
\draw[black] (13+13+2.000000,5.800000) .. controls (13+13+1.000000,7.800000) and (13+13+1.000000,9.800000) .. (13+13+5.000000,9.800000);
\draw[black] (13+13+5.000000,9.800000) .. controls (13+13+6.000000,9.800000) and (13+13+8.000000,7.800000) .. (13+13+7.000000,5.800000);
\draw[blue] (13+13+2,5.8) -- (13+13+2,4);
\node (x) at (13+13+4.5,-1.4) {\huge One-component};
\node (x) at (13+13+4.5,-2.1) {\huge shape $\omega_{k+1}$};

\end{tikzpicture}
\caption{The flip procedure can merge a two-component shape into a one-component shape.}
\label{merge}
\end{figure}

\begin{figure}
  \centering
\subfigure[Two-component initial shape]
  {
    \includegraphics[width=0.31\textwidth,scale=0.35]{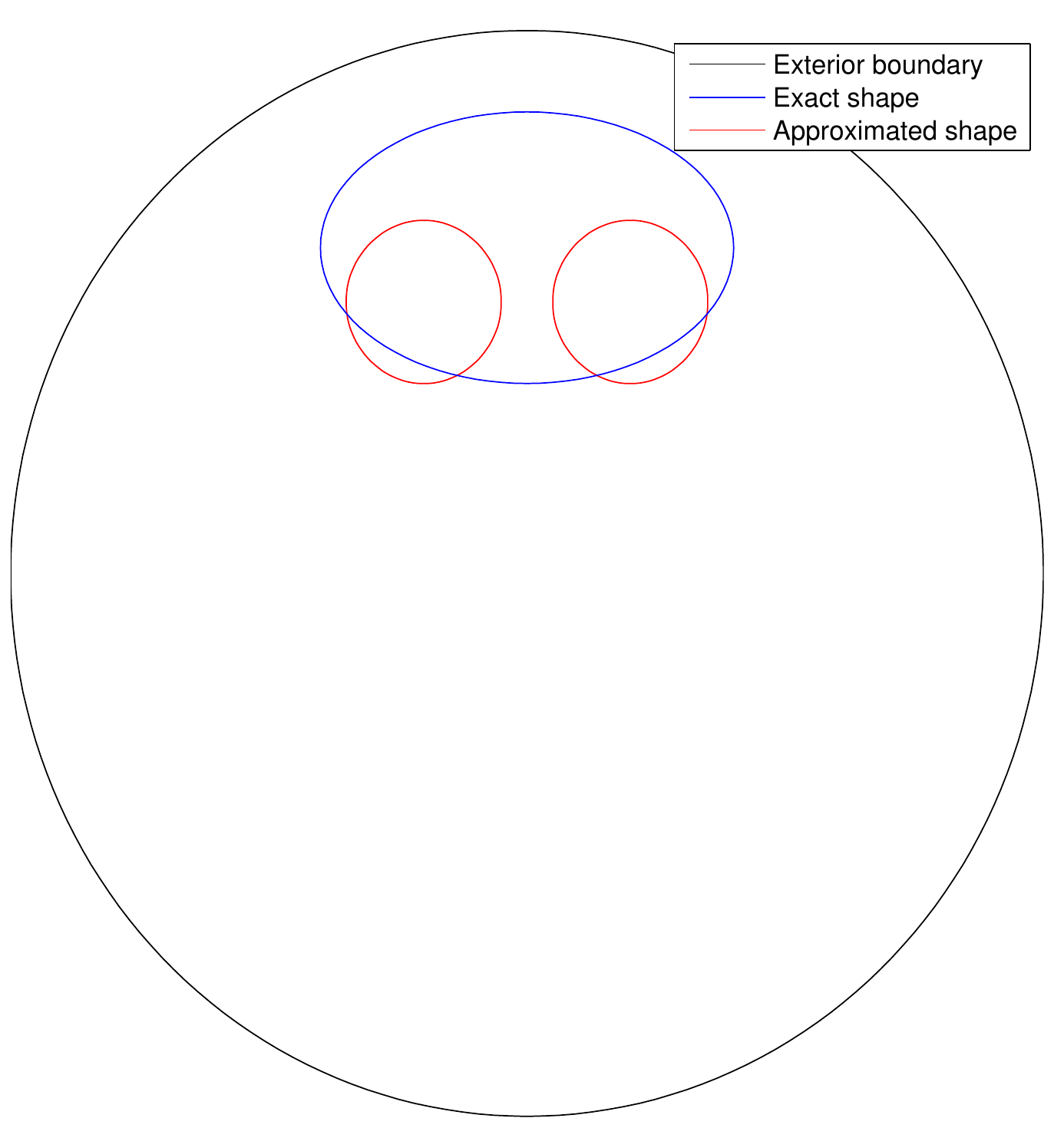}
    \label{simu2shapesstart1}
  }
\quad
\subfigure[Intersecting control polygons]
  {
    \includegraphics[width=0.31\textwidth,scale=0.35]{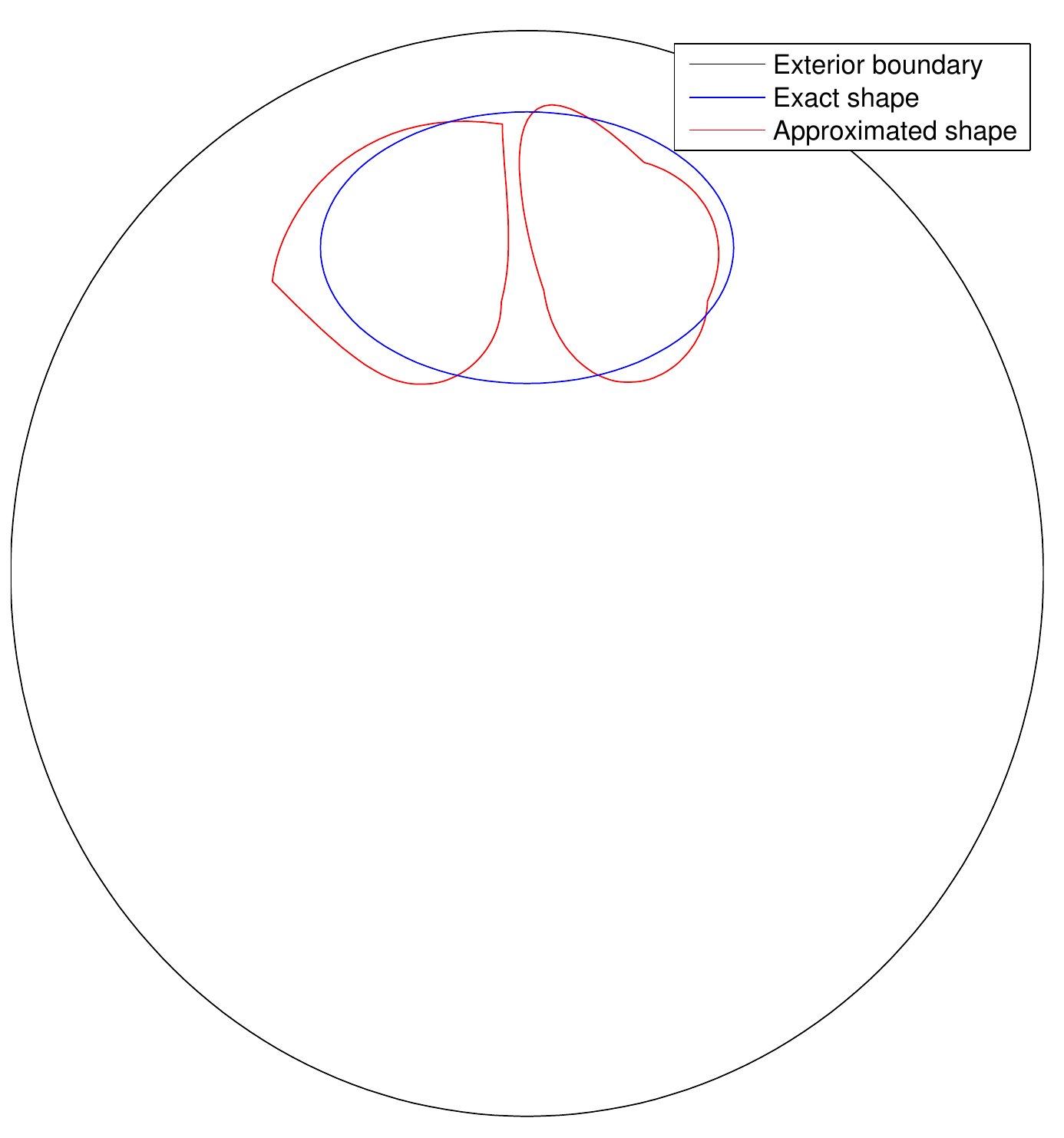}
    \label{simu2shapesstart2}
  }
\subfigure[Flip procedure - The two components have merged]
  {
    \includegraphics[width=0.31\textwidth,scale=0.35]{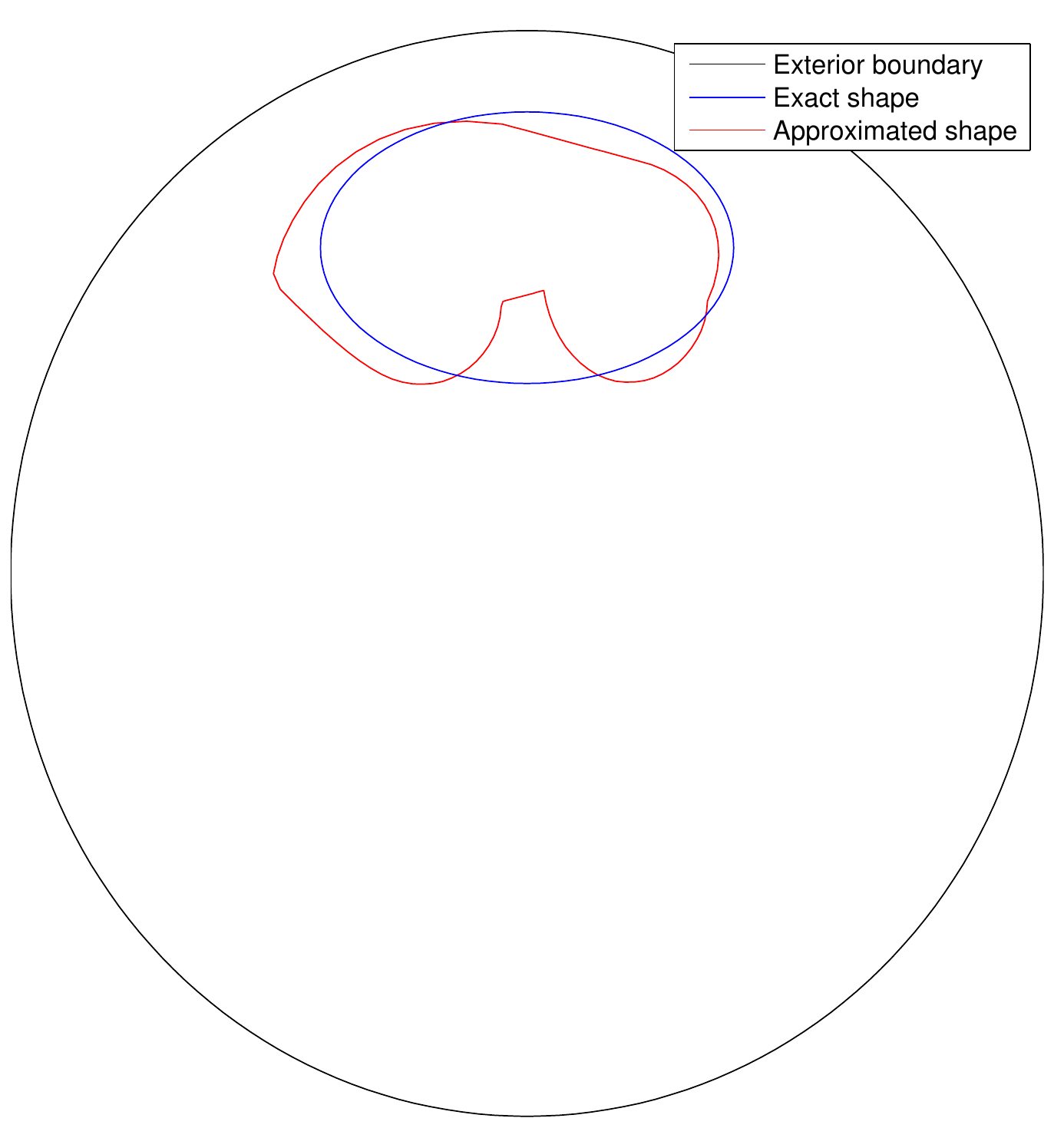}
    \label{simu2shapesstart3}
  }
\quad
\subfigure[Final shape]
  {
    \includegraphics[width=0.31\textwidth,scale=0.35]{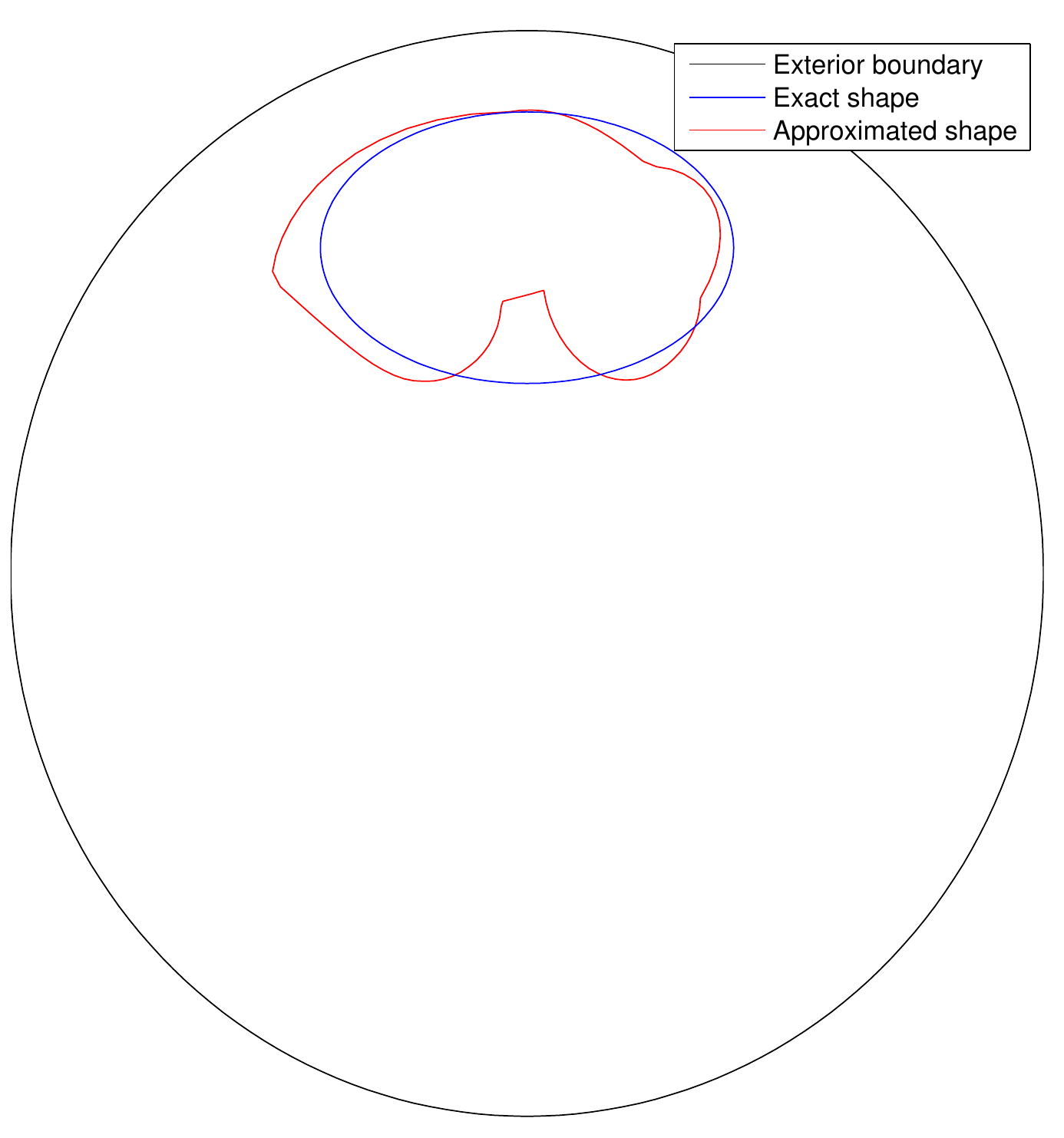}
    \label{simu2shapesstart4}
  }
\caption{Detection of one obstacle starting from a two-component shape}
\label{simuFlipMerge}
\end{figure}

\bibliographystyle{abbrv}

\end{document}